\documentclass[reqno]{amsart}
\usepackage{amssymb, amsmath}
\usepackage{pdfsync}

\theoremstyle{definition}

\theoremstyle{definition} 

\newcommand{\bea}{\begin{eqnarray}}
\newcommand{\eea}{\end{eqnarray}}
\newcommand{\beas}{\begin{eqnarray*}}
\newcommand{\eeas}{\end{eqnarray*}}
\newcommand{\beq}{\begin{equation}}
\newcommand{\eeq}{\end{equation}}

\def\eps{\varepsilon}

\newcommand{\ro}{r}

\newcommand{\cB}{\mathcal B}

\newcommand{\cG}{\mathcal G}

\newcommand{\cN}{\mathcal N}

\newcommand{\cM}{\mathcal{M}}

\newcommand{\R}{\mathbb R}
\newcommand{\RR}{\mathbb R}

\newcommand{\BN}{{\mathbb N}}

\title[On the manifold of closed hypersurfaces]
{On the manifold of closed hypersurfaces in $\R^n$}

\author[Jan Pr{\tiny\"u}{\ss} and Gieri Simonett]{}


\email{jan.pruess@mathematik.uni-halle.de}
\email{gieri.simonett@vanderbilt.edu}

\begin{document}
\maketitle

\centerline{\emph{Dedicated to Jerry Goldstein on the occasion of his 70th anniversary}}
\vspace{1cm}

\centerline{\scshape Jan Pr\"u{\ss}}
\medskip
{\footnotesize
 \centerline{Martin-Luther-Universit\"at Halle-Witten\-berg}
   \centerline{Institut f\"ur Mathematik}
   \centerline{D-06120 Halle, Germany}
} 

\medskip

\centerline{\scshape Gieri Simonett}
\medskip
{\footnotesize
 \centerline{Vanderbilt University}
   \centerline{Department of Mathematics}
   \centerline{Nashville, TN~37240, USA}
}

\bigskip

 \centerline{(Communicated by the associate editor name)}

\begin{abstract}
Several results from differential geometry of hyper-surfaces in $\RR^n$ are derived to form a tool box for the {\em direct mapping method}.
The latter technique has been widely employed to solve problems with moving interfaces, and to study the asymptotics of the induced semiflows.
\end{abstract}


\section{Introduction}\label{intro}
The analysis of problems with moving interfaces has attracted the attention of many researchers
in recent years.
Some of these problems have their origin in mathematical physics, like the Stefan problem,
flows of Newtonian fluids, Hele-Shaw flows, Mullin-Sekerka problems, while others are motivated by problems in differential geometry, like the mean curvature flow, the surface diffusion flow, or the Willmore flow, to mention some prominent examples.

The direct mapping approach to such problems consists in transforming the moving hypersurfaces to a fixed 
reference surface by means of an unknown time-dependent diffeomorphism, which has to be determined as a part of the transformed problem. In the context of the Stefan problem this technique has first been introduced in \cite{Ha81} and is nowadays also called the {\em Hanzawa transform}. The advantage of this approach is that the theory of evolution equations, in particular the theory of maximal regularity, is available for the study of the transformed problems. This way one obtains a local semiflow which, however, does not live in a Banach space as in problems with fixed interfaces, but rather on a manifold which is related to the manifold of hypersurfaces. 
We refer, for instance, to the recent papers
\cite{KPW09, KPW10, PSSS11, PSW11, PSZ10}
by the authors for more details.

To implement this approach one necessarily has to employ results concerning the geometry of hypersurfaces in $\RR^n$, 
and one needs to investigate the structure of the manifold formed by such hypersurfaces.
The main purpose of this paper is to provide a tool box of results 
that are needed for the study of moving interfaces and that are not easily accessible in the literature.
While some of the material presented 
is well-known to researchers specialized in differential geometry
and geometric analysis, we nevertheless believe that the manuscript contains new results
and aspects that are also of interest to specialists.

We investigate the differential geometric properties of embedded hypersurfaces in $n$-dimensional Euclidean space, 
introducing the notion of Weingarten tensor, principal curvatures, mean curvature, tubular neighborhood, surface gradient, surface divergence, and Laplace-Beltrami operator.
The main emphasis lies in deriving representations of these quantities for hypersurfaces $\Gamma=\Gamma_\rho$ that are given as parameterized surfaces in normal direction of a fixed reference surface $\Sigma$
by means of a height function $\rho$.
We derive all of the aforementioned geometric quantities for $\Gamma_\rho$ in terms of $\rho$ and $\Sigma$.
It is also important to study the mapping properties of these quantities in dependence of $\rho$, 
and to derive expressions for their variations.
For instance, we show that 
\begin{equation*}
\kappa^\prime(0)= {\rm tr}\, L_\Sigma^2 +\Delta_\Sigma,
\end{equation*}
where $\kappa=\kappa(\rho)$ denotes the mean curvature of $\Gamma_\rho$,
$L_\Sigma$ the Weingarten tensor of $\Sigma$,
and $\Delta_\Sigma$ the Laplace-Beltrami operator on $\Sigma$.
This is done in Section 3.
In Section~4 we show, among other things, that $C^2$-hypersurfaces can be approximated in a suitable topology 
by smooth (i.e. analytic) hypersurfaces. This leads, in particular, to the existence of parameterizations.
In Section~5 we show that the class of compact embedded hypersurfaces in $\RR^n$
gives rise to a new manifold (whose points are the compact embedded hypersufaces).
Finally, we show that the class $\cM^2(\Omega,r)$ 
of all compact embedded hypersurfaces contained in a bounded domain $\Omega\subset\RR^n$, 
and satisfying a uniform ball condition with radius $r>0$, can be identified with a subspace of $C^2(\bar\Omega)$.
This is important as it allows to derive compactness and embedding properties for 
$\cM^2(\Omega,r)$.    
For further background material in differential geometry we refer to the standard text books in this area, 
e.g.\ to DoCarmo \cite{DoC92} and K\"uhnel~\cite{K02}. 
We also mention \cite{K08} for other aspects on moving hypersurfaces.

\section{Review of some basic differential geometry}
\label{sect-geometry}
We consider a closed embedded hypersurface $\Sigma$ of class $C^k$, $k\geq3$, enclosing a bounded domain $\Omega$
in $\RR^n$. Thus for each point $p\in\Sigma$ there is a ball $B_r(p)\subset\RR^n$ and a diffemorphism $\Phi:B_r(p)\to U\subset\RR^n$ such that
$\Phi(p)=0\in U$ and 
$$\Phi^{-1}(U\cap(\RR^{n-1}\times\{0\}))=B_r(p)\cap\Sigma.$$
We may assume that $\Sigma$ is connected; otherwise we would concentrate on one of its components.
The points of $\Sigma$ are denoted by $p$, and
$\nu_\Sigma=\nu_\Sigma(p)$ means the outer unit normal of $\Sigma$ at $p$. Locally at $p\in\Sigma$ we have the parametrization 
$$p=\phi(\theta):=\Phi^{-1}(\theta,0),$$
where $\theta$ runs through an open parameter set $\Theta\subset\RR^{n-1}$. We denote the tangent vectors generated by this
chart by
\beq\label{tau}
\tau_i=\tau_i(p)=\frac{\partial}{\partial\theta_i}\phi(\theta)=\partial_i\phi, \quad i=1,\ldots,n-1.\eeq
These vectors $\tau_i$ form a basis of the  {\em tangent space} $T_p\Sigma$ of $\Sigma$ at $p$. Note 
that $(\tau_i|\nu_\Sigma)=0$ for all $i$, where $(\cdot |\cdot):=(\cdot |\cdot)_{\R^n}$ 
denotes the Euclidean inner product in $\R^n$.
Similarly, we set $\tau_{ij}=\partial_i\partial_j\phi$, $\tau_{ijk}=\partial_i\partial_j\partial_k \phi$, and so on.
In the sequel we employ Einstein's summation convention, which means that equal lower and upper indices are to be summed,
and $\delta^i_j$ are the entries of the unit matrix $I$. 
For two vectors $a,b\in \R^n$ the tensor product $a\otimes b\in \cB(\R^n)$ is defined by
$[a\otimes b](x)=(b|x)a$ for $x\in\R^n$. If $a$ belongs to the tangent space $T_p\Sigma$,
we may represent $a$ as a linear combination of the basis vectors of $T_p\Sigma$, i.e.
$ a=a^i\tau_i$. The coefficients $a^i$ are called the {\em contravariant components} of $a$. On the other hand, this vector $a$ is also uniquely characterized by its {\em covariant components}, $a_i$ defined by $a_i=(a|\tau_i)$, which means that the covariant components are the coefficients
of the representation of $a$ in the basis  $\{\tau^i\}$ dual to the basis $\{\tau_j\}$, defined by the relations
$(\tau^i|\tau_j)=\delta^i_j$. Similarly,
if $K\in\cB(T_p\Sigma)$ is a tensor we have the representations
$$K=k^{ij}\tau_i\otimes\tau_j=k_{ij}\tau^i\otimes\tau^j=k^i_j\tau_i\otimes\tau^j=k_i^j\tau^i\otimes\tau_j,$$
with e.g.\ $k_{ij}=(\tau_i|K\tau_j)$ and $k^i_j=(\tau^i|K\tau_j)$.


\subsection{The first fundamental form}
Define
\beq\label{gij}g_{ij}=g_{ij}(p)=(\tau_i|\tau_j),\quad i,j=1,\ldots,n-1.\eeq
The matrix $G=[g_{ij}]$ is called {\em the first fundamental form} of $\Sigma$. Note that $G$ is symmetric and also positive definite, since
$$(G\xi|\xi)=g_{ij}\xi^i\xi^j=(\xi^i\tau_i|\xi^j\tau_j)=|\xi^i\tau_i|^2>0, \quad \mbox{ for all } \xi\in\RR^{n-1},\; \xi\neq0.$$
We let $G^{-1}=[g^{ij}]$, hence $g_{ik}g^{kj}=\delta_i^j$, and $g^{il}g_{lj}=\delta^i_j$. The determinant $g:=\det{G}$ is positive.
Let $a$ be a tangent vector. Then $a=a^i\tau_i$ implies
$$a_k=(a|\tau_k)=a^i(\tau_i|\tau_k)=a^ig_{ik}, \quad\mbox { and }\quad a^i=g^{ik}a_k.$$
Thus the fundamental form $G$ allows for the passage from contra- to covariant components of a tangent vector and vice versa.
If $a,b$ are two tangent vectors, then
$$(a|b)=a^ib^j(\tau_i|\tau_j)=g_{ij}a^ib^j=a_jb^j=a^ib_i=g^{ij}a_ib_j=:(a|b)_\Sigma$$
defines an inner product on $T_p\Sigma$ in the canonical way, the {\em Riemannian metric}.
By means of the identity $$(g^{ik}\tau_k|\tau_j)=g^{ik}g_{kj}=\delta^i_j$$
we further see that the dual basis is given by $\tau^i=g^{ik}\tau_k$. We set for the moment $\cG=g^{ij}\tau_i\otimes\tau_j$
and have equivalently
$$\cG= g^{ij}\tau_i\otimes\tau_j=g_{ij}\tau^i\otimes\tau^j= \tau_i\otimes\tau^i=\tau^j\otimes\tau_j.$$
Let $u=u^k\tau_k+(u|\nu_\Sigma)\nu_\Sigma$ be an arbitrary vector in $\RR^n$. Then
$$\cG  u = g^{ij}\tau_i(\tau_j|u)= g^{ij}\tau_i u^kg_{jk}= u^k\tau_k,$$
i.e.\ $\cG$ equals the orthogonal projection $P_\Sigma=I-\nu_\Sigma\otimes\nu_\Sigma$ of $\RR^n$ onto the tangent space $T_p\Sigma$ at $p\in\Sigma$.
Therefore we have the identity
$$ P_\Sigma = I -\nu_\Sigma\otimes\nu_\Sigma = \tau_i\otimes\tau^i=\tau^i\otimes\tau_i.$$
These three properties explain
the meaning of the first fundamental form  $[g_{ij}]$.


\subsection{ The second fundamental form}
Define
\beq\label{lij}l_{ij}=l_{ij}(p)=(\tau_{ij}|\nu_\Sigma), \quad L=[l_{ij}].\eeq
$L$ is called {\em the second fundamental form} of $\Sigma$. Note that $L$ is symmetric, and differentiating the relations
$(\tau_i|\nu_\Sigma)=0$ we derive
\beq\label{lnu} l_{ij}=(\tau_{ij}|\nu_\Sigma)=-(\tau_i|\partial_j\nu_\Sigma)=-(\tau_j|\partial_i\nu_\Sigma).\eeq
The matrix $K$ with entries $l^i_j$, defined by
$$l^i_j=g^{ir}l_{rj},\quad K=G^{-1}L,$$
is  called the {\em shape matrix} of $\Sigma$. The eigenvalues $\kappa_i$ of $K$ are called the {\em principal curvatures}
of $\Sigma$ at $p$,  and the corresponding eigenvectors $\eta_i$ determine the {\em principal curvature directions}. Observe that $K\eta_i=\kappa_i \eta_i$
is equivalent to $L\eta_i=\kappa_i G\eta_i$, hence the relation
$$(L\eta_i|\eta_i) =\kappa_i(G\eta_i|\eta_i)$$ and symmetry of $L$ and $G$ show that the principal curvatures $\kappa_i$ are real. Moreover,
$$\kappa_i(G\eta_i|\eta_j)=(L\eta_i|\eta_j)=(\eta_i|L\eta_j)=\kappa_j(\eta_i|G\eta_j)=\kappa_j(G\eta_i|\eta_j)$$
implies that  principal directions corresponding to different principal curvatures are orthogonal in the Riemannian metric $(\cdot|\cdot)_\Sigma$.
Moreover, the eigenvalues $\kappa_i$ are semi-simple. In fact, if $(K-\kappa_i)x=t\eta_i$ for some $i$ and some $t\in\RR$, then $(L-\kappa_i G)x=tG\eta_i$,
hence
$$ t(G\eta_i|\eta_i)= (Lx-\kappa_iGx|\eta_i)= (x|L\eta_i-\kappa_iG\eta_i)=0,$$
hence $t=0$ since $G$ is positive definite. This shows that $K$ is diagonalizable.

The trace of $K$, i.e.\ the first invariant of $K$, is called the {\em mean curvature} $\kappa$ (times $n-1$) of $\Sigma$ at $p$, i.e. we have
\beq\label{meancurvature}
\kappa_\Sigma={\rm tr}\, K = l^i_i=g^{ij}l_{ij}=\sum_{i=1}^{n-1}\kappa_i.
\eeq
The {\em Gaussian curvature} ${\mathcal K}_\Sigma$ is defined
as the last invariant of $K$, i.e.\
$${\mathcal K}_\Sigma=\det{K} =g^{-1}\det{L} =\Pi_{i=1}^{n-1}\kappa_i.$$
We define the {\em Weingarten tensor} $L_\Sigma$ by means of
\begin{equation}
\label{Weingarten tensor}
L_\Sigma =L_\Sigma(p) 
= l^{ij}\tau_i\otimes\tau_j=l^i_j\tau_i\otimes\tau^j=l_i^j\tau^i\otimes\tau_j=l_{ij}\tau^i\otimes\tau^j.
\end{equation}
$L_\Sigma$ is symmetric with respect to the inner product $(\cdot|\cdot)$ in $\R^n$.
We note that $L_\Sigma\in\cB(\R^n)$ leaves the tangent space $T_p\Sigma$ invariant,
and moreover, that $L_\Sigma \nu_\Sigma=0$.
This shows that $L_\Sigma$ enjoys the direct decomposition
\begin{equation*}
L_\Sigma=L_{\Sigma}|_{T_p\Sigma}\oplus 0: T_p\Sigma \oplus T^{\perp}_p\Sigma \to T_p\Sigma \oplus T^{\perp}_p\Sigma.
\end{equation*}
We will in the following not distinguish between $L_\Sigma$ and its restriction 
$L_{\Sigma}|_{T_p\Sigma}$ to $T_p\Sigma$. 
Observe that
$${\rm tr}L_\Sigma= l^{ij}(\tau_i|\tau_j)=l^{ij}g_{ij}=\kappa_\Sigma,$$ 
and the eigenvalues of $L_\Sigma$  in $T_p\Sigma$ are the principal curvatures since
$$L_\Sigma\eta_k= l^{ij}\tau_i(\tau_j|\eta_k)=l^{ij}\tau_i g_{jr}\eta_k^r=l^i_r\eta_k^r\tau_i=\kappa_k\eta_k^i\tau_i=\kappa_k\eta_k.$$
The remaining eigenvalue of $L_\Sigma$ in $\R^n$ is $0$ with eigenvector $\nu_\Sigma$.


\subsection{The third fundamental form}
To obtain another property of the shape operator $K$ we differentiate the identity $|\nu_\Sigma|^2=1$ to the result $(\partial_i\nu_\Sigma|\nu_\Sigma)=0$.
This means that $\partial_i\nu_\Sigma$ belongs to the tangent space, hence $\partial_i\nu_\Sigma=\gamma^k_i\tau_k$ for some numbers $\gamma_i^k$. Taking the inner
product with $\tau_j$ we get
$$\gamma^k_ig_{kj}=\gamma^k_i(\tau_k|\tau_j)=(\partial_i\nu_\Sigma|\tau_j)=-(\tau_{ij}|\nu_\Sigma)=-l_{ij},$$
hence
$$\gamma^r_i=\gamma^k_ig_{kj}g^{jr}= -l_{ij}g^{jr}=-g^{rj}l_{ji}=-l^r_i,$$
where we used symmetry of $L$ and $G$. Therefore we have
\beq\label{shape}
\partial_i\nu_\Sigma= - l^r_i\tau_r=-l_{ij}\tau^j=-L_\Sigma\tau_i,\quad i=1,\ldots,n-1,
\eeq
the so-called {\em Weingarten relations}.
Furthermore,
$$ 0=\partial_i(\nu_\Sigma|\partial_j\nu_\Sigma)=(\partial_i\nu_\Sigma|\partial_j\nu_\Sigma)
+(\nu_\Sigma|\partial_i\partial_j\nu_\Sigma)$$
implies
\beq\label{secondnu}
-(\partial_i\partial_j\nu_\Sigma|\nu_\Sigma)=(\partial_i\nu_\Sigma|\partial_j\nu_\Sigma)=l^r_il^s_j(\tau_r|\tau_s)
=l^r_ig_{rs}l^s_j=l_{is}g^{sr}l_{rj}=l_i^rl_{rj},
\eeq
which are the entries of the matrix $LG^{-1}L$, i.e.\ the covariant components
of $L_\Sigma^2$. This is the so-called {\em third fundamental form} of $\Sigma$.
In particular this implies the relation
\begin{equation}
{\rm tr}\, L_\Sigma^2= (L_\Sigma\tau^i|L_\Sigma\tau_i)
=- g^{ij}(\partial_i\partial_j\nu_\Sigma|\nu_\Sigma),
\end{equation}
which will be useful later on. Moreover, we deduce from \eqref{secondnu}
\begin{equation}
\label{total curvature}
{\rm tr}\, L_\Sigma^2= (L_\Sigma\tau^i|L_\Sigma\tau_i)
=g^{ij}l^r_i l_{rj}=l^r_i l^i_r=\sum_{i=1}^{n-1} \kappa ^2_i. 
\end{equation}

\subsection{The Christoffel symbols}
The {\em Christoffel symbols} are defined according to
\beq\label{chrisdef}
 \Lambda_{ij|k}=(\tau_{ij}|\tau_k),\quad \Lambda_{ij}^k=g^{kr}\Lambda_{ij|r}.
\eeq
Their importance stems from the representation of $\tau_{ij}$ in the basis $\{\tau_k,\nu_\Sigma\}$ of $\RR^n$ via
\beq\label{chrisprop}
\tau_{ij}=\Lambda_{ij}^k\tau_k +l_{ij}\nu_\Sigma,
\eeq
which follows from $(\nu_\Sigma|\tau_k)=0$ and
$$\Lambda_{ij|k}=(\tau_{ij}|\tau_k)=(\Lambda^r_{ij}\tau_r|\tau_k)=\Lambda^r_{ij}g_{rk}.$$
To express the Christoffel symbols in terms of the fundamental form $G$ we use the identities
\begin{align*}
\partial_kg_{ij} = \partial_k(\tau_i|\tau_j)= (\tau_{ik}|\tau_j)+(\tau_i|\tau_{jk}),\\
\partial_ig_{kj} = \partial_i(\tau_k|\tau_j)= (\tau_{ik}|\tau_j)+(\tau_k|\tau_{ij}),\\
\partial_jg_{ik} = \partial_j(\tau_i|\tau_k)= (\tau_{ij}|\tau_k)+(\tau_i|\tau_{jk}),
\end{align*}
which yield
$$\partial_ig_{jk}+\partial_j g_{ik}-\partial_kg_{ij} = 2(\tau_{ij}|\tau_k),$$
i.e.
\beq\label{chrisbasis}
\Lambda_{ij|k}=\frac{1}{2}[\partial_i g_{jk}+\partial_j g_{ik}-\partial_k g_{ij}].
\eeq
Another important identity follows by differentiation of the relations
$(\tau^j|\tau_k)=\delta^j_k$ and $(\tau^j|\nu_\Sigma)=0$. We have
$$(\partial_i\tau^j|\tau_k)=-(\tau^j|\tau_{ik})= -\Lambda^r_{ik}(\tau^j|\tau_r)
= -\Lambda_{ik}^j,$$
and
$$(\partial_i\tau^j|\nu_\Sigma)=-(\tau^j|\partial_i\nu_\Sigma)=(\tau^j|L_\Sigma\tau_i)=l_i^j,$$
hence
\begin{equation}\label{diffdualbasis}
\partial_i\tau^j= -\Lambda_{ik}^j\tau^k + l_i^j\nu_\Sigma.
\end{equation}
This gives another interpretation of the Christoffel symbols and of
the second fundamental form.


\subsection{The surface gradient}
Let $\rho$ be a scalar field on $\Sigma$. The surface gradient $\nabla_\Sigma\rho$ at $p$ is a vector which belongs to the tangent space of $\Sigma$
at $p$. Thus it can be characterized by its
\begin{itemize}
\item {\em covariant components} $a_i=(\nabla_\Sigma\rho|\tau_i)$, or by its
\item {\em contravariant components} i.e. $\nabla_\Sigma\rho = a^i\tau_i$.
\end{itemize}
The chain rule
$$\partial_i(\rho\circ\phi)=\rho^\prime\partial_i\phi=(\nabla_\Sigma\rho|\tau_i)$$
yields $a_i=\partial_i(\rho\circ\phi)=\partial_i\rho$. This implies
$$a_i=(\nabla_\Sigma\rho|\tau_i)=a^k(\tau_k|\tau_i)=a^kg_{ki},$$
hence
\begin{equation*}
(\nabla_\Sigma\rho)_i=\partial_i(\rho\circ\phi)=\partial_i\rho,
\quad (\nabla_\Sigma\rho)^i=g^{ij}\partial_j\rho,\quad
\nabla_\Sigma\rho= \tau^i\partial_i\rho,
\quad
\nabla_\Sigma\rho=(g^{ij}\partial_j\rho)\tau_i.
\end{equation*}
For a scalar field $\rho$ defined in a neighborhood of $\Sigma$ we therefore have
$$\nabla\rho=(\nabla\rho|\nu_\Sigma)\nu_\Sigma +(\nabla_\Sigma\rho)^i\tau_i,$$
and hence, the surface gradient of $\rho$ is the projection of $\nabla\rho$ onto $T_p\Sigma$, that is,
$$\nabla_\Sigma\rho=P_\Sigma\nabla\rho.$$
For a vector field $f:\Sigma\rightarrow\RR^m$ of class $C^1$ we define similarly
$$\nabla_\Sigma f := g^{ij}\tau_i\otimes \partial_j f=\tau^j\otimes\partial_j f.$$
In particular, this yields for the identity map ${\rm id}_\Sigma$ on $\Sigma$
$$\nabla_\Sigma\,{\rm id}_\Sigma=g^{ij}\tau_i\otimes\partial_j\phi= g^{ij}\tau_i\otimes\tau_j= P_\Sigma,$$
and by the Weingarten relations
$$\nabla_\Sigma\nu_\Sigma= g^{ij}\tau_i\otimes\partial_j\nu_\Sigma= - g^{ij}l_j^r\tau_i\otimes\tau_r=
-l^{ij}\tau_i\otimes\tau_j=-L_\Sigma.$$
For the surface gradient of tangent vectors we have
\begin{eqnarray*}
\nabla_\Sigma \tau_k&=& g^{ij}\tau_i\otimes \partial_j\tau_k=g^{ij}\tau_i\otimes\tau_{jk}
= g^{ij}\tau_i\otimes(\Lambda_{jk}^r\tau_r+l_{jk}\nu_\Sigma)\\
&=& g^{ij}\Lambda_{jk}^r\tau_i\otimes\tau_r + l^i_k\tau_i\otimes\nu_\Sigma
= \Lambda_{kj}^r\tau^j\otimes\tau_r + (L_\Sigma\tau_k)\otimes\nu_\Sigma.
\end{eqnarray*}


\subsection {The surface divergence}
Let $f$ be a tangential vector field on $\Sigma$. 
As before, $f^i=(f|\tau^i)$ denote the contravariant components of $f$,
and $f_i=(f|\tau_i)$ the covariant components, respectively.
The {\em surface divergence} of $f$ is defined by
\beq\label{surfdiv}
{\rm div}_\Sigma\, f = \frac{1}{\sqrt{g}}\partial_i(\sqrt{g}f^i)= \frac{1}{\sqrt{g}}\partial_i(\sqrt{g}g^{ij}f_j).
\eeq
As before, $g:=\det G$ denotes the determinant of $G=[g_{ij}]$.
This definition ensures that partial integration can be carried out as usual:
$$\int_\Sigma (\nabla_\Sigma \rho|f)_\Sigma \,d\sigma= -\int_\Sigma \rho\,{\rm div}_\Sigma f \,d\sigma.$$
Recall that the surface measure in local coordinates is given by $d\sigma=\sqrt{g}d\theta$, which explains the factor $\sqrt{g}$.

In fact, if e.g.\ $\rho$ has support in a chart $\phi(\Theta)$ at $p$ then
\begin{align*}
\int_\Sigma (\nabla_\Sigma \rho|f)_\Sigma \,d\sigma
&= \int_\Theta \partial_i(\rho\circ\phi) [(f^i\circ\phi) \sqrt{g})]\,d\theta\\
&= -\int_\Theta (\rho\circ\phi)\frac{1}{\sqrt{g}}\partial_i [\sqrt{g}(f^i\circ\phi)]\sqrt{g}\,d\theta
= -\int_\Sigma \rho\, {\rm div}_\Sigma f \,d\sigma.
\end{align*}
There is another useful representation of surface divergence, given by
\beq
\label{gendiv}
{\rm div}_\Sigma f = g^{ij}(\tau_j |\partial_if)
=(\tau^i|\partial_if).
\eeq
It comes from
$${\rm div}_\Sigma f = \frac{1}{\sqrt{g}}\partial_i(\sqrt{g}g^{ij}f_j)
=\frac{1}{\sqrt{g}}\partial_i[\sqrt{g}g^{ij}(\tau_j|f)],$$
since
\begin{equation}
\label{div-0}
(\partial_i(\sqrt{g}g^{ij}\tau_j)|\tau_k)=0,\quad k=1,\ldots,n-1.
\end{equation}
Here \eqref{div-0} follows from 
\begin{equation*}
\begin{split}
(\partial_i(\sqrt{g}g^{ij}\tau_j)|\tau_k)
&=\partial_i(\sqrt{g}g^{ij}(\tau_j|\tau_k))-\sqrt{g}g^{ij}(\tau_j|\tau_{ki})
=\partial_k\sqrt{g}-\sqrt{g}g^{ij}(\tau_j|\tau_{ki})\\
&=\partial_k\sqrt{g}-\frac{1}{2}\sqrt{g}g^{ij}\partial_k(\tau_j|\tau_i)
=\frac{1}{2\sqrt{g}}\left(\partial_k g- g g^{ij}\partial_k g_{ij}\right) \\
\end{split}
\end{equation*}
and the well-known relation
\begin{equation*}
\begin{split}
\partial_k g=\partial_k \det G\,
&=\sum_{j=1}^n\det\,[g_{\bullet 1},\cdots,\partial_k g_{\bullet j},\cdots g_{\bullet n}]\\
&=(\det G)\sum_{j=1}^n\det\big(G^{-1}[g_{\bullet 1},\cdots,\partial_k g_{\bullet j},\cdots g_{\bullet n}]\big)\\
&=g\,{\rm tr}\,[G^{-1}\partial_k G]=g g^{ij}\partial_k g_{ij}, 
\end{split}
\end{equation*}
where $G=[g_{ij}]=[g_{\bullet 1},\cdots,g_{\bullet n}]$, with $g_{\bullet j}$ the $j$-th column of $G$.

Equation (\ref{gendiv}) can be be used as a definition of surface divergence for general,
not necessarily tangential vector fields $f$.
For example, consider $f=\nu_\Sigma$; then $\partial_i\nu_\Sigma=-l_i^k\tau_k$ by the Weingarten relations, hence we obtain
\begin{equation*}
{\rm div}_\Sigma\nu_\Sigma = g^{ij}(\tau_j|\partial_i\nu_\Sigma)
=-g^{ij}l_{ij}=-\kappa_\Sigma.
\end{equation*}
This way we have derived the important relation
\beq\label{mcdiv}
\kappa_\Sigma=-{\rm div}_\Sigma\nu_\Sigma.\eeq
Note that the surface divergence theorem only holds for tangential vector fields!

Another representation of the surface divergence of a general vector field $f$ is given by
$${\rm div}_\Sigma f =(\tau^i|\partial_if)= {\rm tr}[\tau^i\otimes\partial_i f]
={\rm tr}\, \nabla_\Sigma f.$$
Finally, we compute
$${\rm div}_\Sigma \tau_k = g^{ij}(\tau_j|\tau_{ki}) = g^{ij}\Lambda_{ki|j}= \Lambda_{ik}^i.$$


\subsection{The Laplace-Beltrami operator}
The {\em Laplace-Beltrami} operator on $\Sigma$ is defined for scalar fields by means of
$$ \Delta_\Sigma \rho = {\rm div}_\Sigma\nabla_\Sigma \rho,$$
which in local coordinates reads
$$\Delta_\Sigma \rho = \frac{1}{\sqrt{g}}\partial_i[\sqrt{g}g^{ij}\partial_j\rho].$$
Another representation of $\Delta_\Sigma$ is given by
\begin{equation}
\label{LB-Christoffel}
\Delta_\Sigma\rho =g^{ij}\partial_i\partial_j\rho-g^{ij}\Lambda_{ij}^k\partial_k\rho.
\end{equation}
In order to see this we note that \eqref{div-0} implies
\begin{equation*}
0=(\partial_i(\sqrt g g^{ir}\tau_r)|\tau_k)=\partial_i(\sqrt{g}g^{ir})g_{rk}
+\sqrt{g} g^{ir}(\tau_{ir}|\tau_k)
\end{equation*}
and hence
\begin{equation*}
\partial_i (\sqrt{g}g^{ij})=-\sqrt{g}g^{ir} g^{jk}(\tau_{ir}|\tau_k)=-\sqrt{g}g^{ir}\Lambda^j_{ir}.
\end{equation*}
Since at each point $p\in\Sigma$ we may choose a chart such that $g_{ij}=\delta_{ij}$
and $\Lambda_{ij}^k=0$ at $p$, we see from this representation
that the Laplace-Beltrami operator is equivalent to the Laplacian at the point $p$;
see subsection 8 below.

To obtain another representation of $\Delta_\Sigma$, for a scalar $C^2$-function we compute
$$\nabla_\Sigma^2\rho=\nabla_\Sigma(\tau^j\partial_j\rho)=
\tau^i\otimes\partial_i(\tau^j\partial_j\rho).$$
This yields with 
\eqref{diffdualbasis}
\begin{align*}
\nabla_\Sigma^2\rho=& (\partial_i\partial_j\rho)\tau^i\otimes\tau^j +
(\partial_j\rho)\tau^i\otimes\partial_i\tau^j\\
=& (\partial_i\partial_k\rho-\Lambda^j_{ik}\partial_j\rho)\tau^i\otimes\tau^k
+ (L_\Sigma\nabla_\Sigma\rho)\otimes\nu_\Sigma.
\end{align*}
Taking traces gives
$$\Delta_\Sigma \rho={\rm tr}\, \nabla_\Sigma^2 \rho.$$
Similarly, the Laplace-Beltrami operator applies to general vector fields $f$ according to
$$\Delta_\Sigma f = g^{ij}(\partial_i\partial_j f -\Lambda_{ij}^r\partial_r f).$$
For example, this yields for the identity map ${\rm id}_\Sigma $ on $\Sigma$ 
$$\Delta_\Sigma\,{\rm id}_\Sigma  = g^{ij}(\partial_i\partial_j \phi-\Lambda_{ij}^r\partial_r\phi)= g^{ij}(\tau_{ij}-\Lambda_{ij}^r\tau_r),$$
and hence by \eqref{chrisprop}
$$\Delta_\Sigma\, {\rm id}_\Sigma= g^{ij}l_{ij}\nu_\Sigma = \kappa_\Sigma\nu_\Sigma.$$
Finally, we prove the important formula
\begin{align} \label{gradcurv}
\Delta_\Sigma\nu_\Sigma=-\nabla_\Sigma\kappa_\Sigma -[{\rm tr}\, L_\Sigma^2]\nu_\Sigma.
\end{align}
In fact, we have from \eqref{secondnu}
\begin{equation*}
(\Delta_\Sigma\nu_\Sigma|\nu_\Sigma)
=g^{ij}(\partial_{ij}\nu_\Sigma-\Lambda^r_{ij}\partial_r\nu_\Sigma|\nu_\Sigma)
= g^{ij}(\partial_{ij}\nu_\Sigma|\nu_\Sigma)=-{\rm tr}\, L_\Sigma^2.
\end{equation*}
Next observe that
\begin{align*}
&(\partial_k\partial_j\nu_\Sigma|\tau_i)-(\partial_i\partial_j\nu_\Sigma|\tau_k)
=\partial_k(\partial_j\nu_\Sigma|\tau_i)-\partial_i(\partial_j\nu_\Sigma |\tau_k)\\
&=-\partial_k(\nu_\Sigma|\tau_{ij})+\partial_i(\nu_\Sigma|\tau_{kj})= \partial_k(\partial_i\nu_\Sigma|\tau_j)
-\partial_i(\partial_k\nu_\Sigma|\tau_j)\\
&=(\partial_i\nu_\Sigma|\tau_{kj})-(\partial_k\nu_\Sigma|\tau_{ij})= \Lambda^r_{kj}(\partial_i\nu_\Sigma|\tau_r)
-\Lambda^r_{ij}(\partial_k\nu_\Sigma|\tau_r)\\
&=\Lambda^r_{kj}(\partial_r\nu_\Sigma|\tau_i)-\Lambda^r_{ij}(\partial_r\nu_\Sigma|\tau_k),
\end{align*}
hence
$$(\partial_k\partial_j\nu_\Sigma-\Lambda_{kj}^r\partial_r\nu_\Sigma|\tau_i)=
(\partial_i\partial_j\nu_\Sigma-\Lambda_{ij}^r\partial_r\nu_\Sigma|\tau_k).$$
This implies
$$(\Delta_\Sigma\nu_\Sigma|\tau_i)= g^{jk}(\partial_k\partial_j\nu_\Sigma-\Lambda_{kj}^r\partial_r\nu_\Sigma|\tau_i)
=(\partial_i\partial_j\nu_\Sigma-\Lambda_{ij}^r\partial_r\nu_\Sigma|\tau^j).$$
On the other hand,
\begin{align*}
-(\nabla_\Sigma\kappa_\Sigma|\tau_i)&= -\partial_i\kappa_\Sigma =
\partial_i(\partial_j\nu_\Sigma|\tau^j)\\
&= (\partial_i\partial_j\nu_\Sigma|\tau^j)+ (\partial_r\nu_\Sigma|\partial_i\tau^r)\\
&=(\partial_i\partial_j\nu_\Sigma - \Lambda^r_{ij}\partial_r\nu_\Sigma|\tau^j).
\end{align*}
This proves formula \eqref{gradcurv}.


\subsection{The case of a graph over $\RR^{n-1}$}
Suppose that $\Sigma$ is a graph over $\RR^{n-1}$, i.e.\ there is a function $h\in C^2(\RR^{n-1})$ such that the hypersurface $\Sigma$
is given by the chart $\phi(x)=[x^T,h(x)]^T$, $x\in\RR^{n-1}$. Then the tangent vectors are given by $\tau_i=[e_i^T,\partial_i h]^T$,
where $\{e_i\}$ denotes the standard basis in $\RR^{n-1}$. The (upward pointing) normal $\nu_\Sigma$ is given by
$$\nu_\Sigma(x)= \beta(x)[-\nabla_xh(x)^T,1]^T,\quad \beta(x)=1/\sqrt{1+|\nabla_x h|^2}  .$$
The first fundamental form becomes
$$g_{ij}= \delta_{ij} + \partial_i h\partial_j h,$$
hence
$$g^{ij}= \delta^{ij} - \beta^2\partial_ih\partial_j h.$$
This yields
$$\tau^i =[[e_i -\beta^2\partial_ih\nabla_xh]^T,\beta^2\partial_ih]^{\sf T},$$
and with $\tau_{ij}=[0,\partial_i\partial_j h]^T$
$$l_{ij}= (\tau_{ij}|\nu_\Sigma)= \beta \partial_{i}\partial_jh,$$
and therefore
$$ \kappa_\Sigma= g^{ij}l_{ij}= \beta[\Delta_xh- \beta^2(\nabla^2_xh\nabla_xh|\nabla_xh)]=
{\rm div}_x\frac{\nabla_xh}{\sqrt{1+|\nabla_xh|^2}}.$$
The Christoffel symbols  in this case are given by
$$ \Lambda_{ij|k}=\partial_i\partial_jh\partial_kh,\quad \Lambda_{ij}^k
= \beta^2\partial_i\partial_jh\partial_kh.$$
Suppose that $\RR^{n-1}\times\{0\}$ is the tangent plane at $\phi(0)=0\in\Sigma$. Then $h(0)=\nabla_xh(0)=0$, hence at this
point we have
$g_{ij}=\delta_{ij}$, $\tau_i=[e_i^T,0]^T$, $\nu_\Sigma=[0,1]^T$, $\beta=1$, and $l_{ij}=\partial_i\partial_jh$. Thus the curvatures
are the eigenvalues of $\nabla^2_xh$, the mean curvature is $\kappa_\Sigma=\Delta_xh$, and $\Lambda_{ij}^k=0$.

To obtain a representation of the surface gradient, let $\rho:\Sigma\to\RR$, then
$$\nabla_\Sigma\rho = \tau^j\partial_j\rho 
= [[\nabla_x\rho-\beta^2(\nabla_x\rho|\nabla_xh)\nabla_xh]^T,\beta^2 (\nabla_x\rho|\nabla_xh)]^{\sf T}.$$
Similarly, for $f=(\bar f, f^n):\Sigma\to\R^{n-1}\times\R $ we obtain
$${\rm div}_\Sigma f = (\tau^i|\partial_i f)= {\rm div}_x\bar f
+\beta^2(\nabla_x h|\nabla_xf^n-(\nabla_xh\cdot\nabla_x)\bar f),$$
and for the Laplace-Beltrami
$$\Delta_\Sigma\rho= \Delta_x\rho -\beta^2(\nabla_x^2\rho\nabla_xh|\nabla_xh)
-\beta^2[\Delta_xh -\beta^2(\nabla_x^2h\nabla_xh|\nabla_xh)](\nabla_xh|\nabla_x\rho).$$

\medskip

\section{Parameterized hypersurfaces}
We consider now a hypersurface $\Gamma=\Gamma_\rho$ which is parameterized over a fixed hypersurface $\Sigma$ according to
\beq\label{gammarho}
q=\psi_\rho(p)= p+\rho(p)\nu_\Sigma(p),\quad p\in\Sigma,
\eeq
where as before $\nu_\Sigma=\nu_\Sigma(p)$ denotes the outer normal of $\Sigma$ at $p\in\Sigma$. We want to derive the basic geometric quantities
of $\Gamma$ in terms of $\rho$ and those of $\Sigma$. In the sequel we assume that $\rho$ is of class $C^1$ and small enough. A precise bound on $\rho$ will be given below.


\subsection{The fundamental form}
Differentiating (\ref{gammarho}) we obtain
$$\tau_i^\Gamma=\partial_i\psi_\rho = \tau_i +\rho\partial_i\nu_\Sigma + (\partial_i\rho)\nu_\Sigma,$$
hence the Weingarten relations imply
\beq\label{taugamma}
\tau_i^\Gamma=\tau_i -\rho l^r_i\tau_r +\nu_\Sigma\partial_i\rho= (I-\rho L_\Sigma)\tau_i +\nu_\Sigma\partial_i\rho.\eeq
Therefore we may compute the fundamental form $G^\Gamma=[g^\Gamma_{ij}]$ of $\Gamma$.
\begin{align*}
g^\Gamma_{ij}&=(\tau^\Gamma_i|\tau^\Gamma_j)=(\tau_i|\tau_j)-(\tau_i|\rho l^r_j\tau_r)+(\tau_i|\nu_\Sigma)\partial_j\rho\\
&-(\rho l^r_i\tau_r|\tau_j)+(\rho l^r_i\tau_r|\rho l^s_j\tau_s)-(\rho l^r_i\tau_r|\nu_\Sigma)\partial_j\rho\\
&+(\nu_\Sigma|\tau_j)\partial_i\rho-(\nu_\Sigma|l^s_j\tau_s)\partial_i\rho+|\nu_\Sigma|^2\partial_i\rho\partial_j\rho.
\end{align*}
Since $|\nu_\Sigma|^2=1$ and $(\nu_\Sigma|\tau_i)=0$, this yields
\beq\label{gijgamma}
g_{ij}^\Gamma=g_{ij}-2\rho l_{ij} +\rho^2l^r_il_{rj} +\partial_i\rho\partial_j\rho.
\eeq
Let
$$g_{ij}(\rho)=g_{ij}-2\rho l_{ij} +\rho^2l^r_il_{rj};$$
then we may write
$$g^\Gamma_{ij}=g_{ik}(\rho)[\delta^k_j+ g^{kr}(\rho)\partial_r\rho\partial_j\rho],$$
hence
$$ G^\Gamma= G(\rho)(I+G^{-1}(\rho)\partial\rho\otimes\partial\rho)$$
where $\partial\rho=[\partial_1\rho,\dots,\partial_{n-1}\rho]^T$.
Next we may factor $G(\rho)$ according to
$$G(\rho)=G(I-2\rho G^{-1}L+ \rho^2 (G^{-1}L)^2) = G(I-\rho K)^2.$$
Since for any two vectors $a,b\in\RR^n$ we have
$$\det(I+a\otimes b ) =1+(a|b),$$
we obtain
\beq\label{detggamma}
g^\Gamma:=\det{G^\Gamma} = g[\det(I-\rho K)]^2(1+(G^{-1}(\rho)\partial\rho|\partial\rho))=g\alpha^2(\rho)/\beta^2(\rho),\eeq
where
$$\alpha(\rho)=\det(I-\rho K)=\Pi_{i=1}^{n-1}(1-\rho\kappa_i),$$
and
$$\beta(\rho)=1/\sqrt{1+(G^{-1}(\rho)\partial\rho|\partial\rho)}.$$
This yields for the surface measure $d\gamma$ on $\Gamma_\rho$
\beq\label{smgamma}
d\gamma=\sqrt{g^\Gamma}d\theta = \frac{\alpha(\rho)}{\beta(\rho)}\sqrt{g}\,d\theta
=\frac{\alpha(\rho)}{\beta(\rho)}\,d\sigma,
\eeq
hence
$${\rm mes}\,\Gamma_\rho =\int_{\Gamma_\rho} d\gamma=\int_\Sigma \frac{\alpha(\rho)}{\beta(\rho)}\,d\sigma.$$
Since
$$(I+a\otimes b)^{-1}= I - \frac{a\otimes b}{1+(a|b)},$$
we obtain for $[G^\Gamma]^{-1}$ the identity
\beq\label{ggammainv}
[G^\Gamma]^{-1}=[I-\beta^{2}(\rho)G^{-1}(\rho)\partial\rho\otimes\partial\rho)]G^{-1}(\rho),
\eeq
and
$$G^{-1}(\rho)=(I-\rho K)^{-2}G^{-1}.$$
All of this makes sense only for functions $\rho$ such that $I-\rho K$ is invertible, i.e.\ $\alpha(\rho)$ should not vanish.
Thus the precise bound for $\rho$ is determined by the principle curvatures of $\Sigma$, 
and we assume here and in the sequel that
\beq
\label{boundrho}
|\rho|_\infty< \frac{1}{\max\{|\kappa_i(p)|:\,i=1,\ldots,n-1,\, p\in\Sigma\}}=:\rho_0.\eeq


\subsection{The normal at $\Gamma$}
We next compute the unit outer normal at $\Gamma$.
For this purpose we set
$$\nu^\Gamma=\beta(\rho)(\nu_\Sigma-a(\rho)),$$
where $\beta$ is a scalar and $a(\rho)\in T_p\Sigma$.
Then $\beta(\rho)=(1+|a(\rho)|^2)^{-1/2}$ and
$$0=(\nu^\Gamma|\tau_i^\Gamma)/\beta(\rho)= (\nu_\Sigma-a|\tau_i-\rho L_\Sigma\tau_i+\nu_\Sigma\partial_i\rho),$$
which yields
$$0= \partial_i\rho -(a(\rho)|(I-\rho L_\Sigma)\tau_i)=\partial_i\rho -((I-\rho L_\Sigma)a(\rho)|\tau_i),$$
by symmetry of $L_\Sigma$. But this implies $(I-\rho L_\Sigma)a(\rho)=\nabla_\Sigma\rho$, i.e.\ we have
\begin{equation}\label{normalGamma}
\nu^\Gamma= \beta(\rho)(\nu_\Sigma-M_0(\rho)\nabla_\Sigma\rho),
\end{equation}
with
$$M_0(\rho)=(I-\rho L_\Sigma)^{-1},\quad
\beta(\rho)=(1+|M_0(\rho)\nabla_\Sigma\rho|^2)^{-1/2}.$$
As remarked in subsection 2.2 we do not distinguish between $L_\Sigma\in\cB(\R^n)$ and its restriction to $T_p\Sigma$.
With this identification, and by the fact that
 $(I-\rho L_\Sigma )=I$ on $T^\perp_p\Sigma $,
we have
\begin{equation*}
(I-\rho  L_\Sigma)(p)\in{\rm Isom}\,(\R^n,\R^n)\cap {\rm Isom}\,(T_p\Sigma,T_p\Sigma),
\end{equation*}
provided $\rho$ satisfies \eqref{boundrho}.
As before, $\rho L_\Sigma$ is short form for $\rho(p)L_\Sigma(p)$.
Hence, we have
\begin{equation*}
M_0(\rho)(p)\in {\rm Isom}(\R^n,\R^n)\cap {\rm Isom}(T_p\Sigma,T_p\Sigma).
\end{equation*}
Note that $\beta(\rho)$ coincides with $\beta(\rho)$ as defined in the previous subsection.
By means of $a(\rho)$, $\beta(\rho)$ and $M_0(\rho)$ this leads to another 
representation of $G^\Gamma$ and $G_\Gamma^{-1}$, namely
$$G^\Gamma = G_\Sigma(I-\rho L_\Sigma)[I+a(\rho)\otimes a(\rho)](I-\rho L_\Sigma),$$
and
$$G_\Gamma^{-1} = M_0(\rho)[I-\beta^2(\rho)a(\rho)\otimes a(\rho)]M_0(\rho)G_\Sigma^{-1}.$$


\subsection{The surface gradient and the surface divergence on $\Gamma$}
It is  of importance to have a representation for the surface gradient on $\Gamma$ in
terms of $\Sigma$. For this purpose recall that
$$ P_\Gamma =I-\nu^\Gamma\otimes\nu^\Gamma =g^{ij}_\Gamma \tau_i^\Gamma\otimes \tau_j^\Gamma,$$
where $\nu^\Gamma =\beta(\rho)(\nu_\Sigma -M_0(\rho)\nabla_\Sigma \rho)$, and
$$\tau_i^\Gamma= (I-\rho L_\Sigma)\tau_i^\Sigma + (\nabla_\Sigma\rho)_i\nu_\Sigma.$$
By virtue of $L_\Sigma\nu_\Sigma =0$, the latter implies
$$M_0(\rho)\tau^\Gamma_i = \tau_i^\Sigma +(\partial_i\rho)\nu_\Sigma,$$
hence \begin{equation}\label{B1}
P_\Sigma M_0(\rho) \tau^\Gamma_i =\tau_i^\Sigma.
\end{equation}
On the other hand, we have
$$P_\Gamma M_0(\rho)\tau^r_\Sigma = g^{ij}_\Gamma \tau^\Gamma_i\otimes
\tau^\Gamma_j M_0(\rho)\tau^r_\Sigma  = \tau^j_\Gamma (\tau_j^\Gamma|M_0(\rho)\tau^r_\Sigma),$$
hence
\begin{equation}\label{B2}P_\Gamma M_0(\rho)\tau^r_\Sigma = \tau^j_\Gamma(M_0(\rho)\tau_j^\Gamma|\tau^r_\Sigma)=
\tau^r_\Gamma.\end{equation}
\eqref{B1} and \eqref{B2} allow for an easy change between the bases of $T_p\Sigma$ and $T_q\Gamma$,
where $q=\psi_\rho(p)=p+\rho(p)\nu_\Sigma(p)$.

\eqref{B2} implies for a scalar function $\varphi$ on $\Gamma$
$$\nabla_\Gamma\varphi =\tau^r_\Gamma \partial_r\varphi = P_\Gamma M_0(\rho)
\tau_\Sigma^r\partial_r\varphi_* = P_\Gamma M_0(\rho)\nabla_\Sigma \varphi_*,\quad \varphi_*=\varphi\circ\psi_\rho$$
which leads to the identity
$$ \nabla_\Gamma \varphi =P_\Gamma M_0(\rho) \nabla_\Sigma \varphi_*.$$
Similarly, if $f$ denotes a vector field on $\Gamma$ then
$$\nabla_\Gamma f = P_\Gamma M_0(\rho)\nabla_\Sigma f_*,$$
and so
$${\rm div}_\Gamma f = (\tau^r_\Gamma|\partial_r f)=
(P_\Gamma M_0(\rho) \tau^r_\Sigma|\partial_r f) =
{\rm tr}\, [P_\Gamma M_0(\rho)\nabla_\Sigma f_*].$$
As a consequence, we obtain for the Laplace-Beltrami operator on $\Gamma$
$$\Delta_\Gamma \varphi = {\rm tr}\, [P_\Gamma M_0(\rho)\nabla_\Sigma (P_\Gamma
M_0(\rho)\nabla_\Sigma\varphi_*)],$$
which can be written as
$$\Delta_\Gamma\varphi= M_0(\rho)P_\Gamma M_0(\rho):\nabla_\Sigma^2\varphi_* 
+ (b(\rho,\nabla_\Sigma\rho,\nabla^2_\Sigma\rho)|\nabla_\Sigma\varphi_*),$$
with $b=\partial_i(M_0P_\Gamma)M_0P_\Gamma\tau^i_\Sigma$.
One should note that the structure of the Laplace-Beltrami operator on $\Gamma$ in local coordinates is
\begin{equation*}
 \Delta_\Gamma \varphi 
= a^{ij}(\rho,\partial\rho) \partial_i\partial_j\varphi_*
+ b^k(\rho,\partial\rho,\partial^2\rho)\partial_k\varphi_*
\end{equation*}
with
\begin{equation*}
a^{ij}(\rho,\partial\rho)
=(P_\Gamma M_0(\rho)\tau^i_\Sigma |P_\Gamma M_0(\rho)\tau^j_\Sigma)
=(\tau^i_\Gamma |\tau^j_\Gamma)=g^{ij}_\Gamma
\end{equation*}
and
\begin{equation*}
b^k(\rho,\partial\rho,\partial^2\rho))
=(\partial_i(M_0(\rho)P_\Gamma)P_\Gamma M_0(\rho)\tau^i_\Sigma|\tau^k_\Sigma)
=(\tau^i_\Gamma|\partial_i \tau^k_\Gamma)
=-g^{ij}_\Gamma \Lambda^k_{\Gamma ij}. 
\end{equation*}
This shows that $-\Delta_\Gamma$ is strongly elliptic on the reference manifold $\Sigma$
as long as $|\rho|_\infty<\rho_0$. 


\subsection{Normal variations.}
For $\rho, h\in C(\Sigma)$ sufficiently smooth  and a mapping $M(\rho)$ we define
\begin{equation*}
M^\prime (0) h:=\frac{d}{d\varepsilon}M(\rho+\varepsilon h)\Big|_{\varepsilon=0}.
\end{equation*}
First we have
$$M_0^\prime(\rho)= M_0(\rho)L_\Sigma M_0(\rho),\quad  M_0^\prime(0)=L_\Sigma,$$
as $M_0(0)=I$. Next
$$\beta^\prime(\rho)h=-\beta(\rho)^3
\big(M_0(\rho)\nabla_\Sigma\rho \big|M_0^\prime(\rho)h\nabla_\Sigma \rho +M_0(\rho)\nabla_\Sigma h\big),$$
which yields $\beta^\prime(0)=0$, as $\beta(0)=1$.
From this we get for the normal $$\nu(\rho)=\nu^\Gamma=\beta(\rho)(\nu_\Sigma-M_0(\rho)\nabla_\Sigma\rho)$$ the relation
$$\nu^\prime(\rho)h=\beta^\prime(\rho)h(\nu_\Sigma-M_0(\rho)\nabla_\Sigma\rho)-\beta(\rho)(M_0^\prime(\rho)h\nabla_\Sigma\rho +M_0(\rho)\nabla_\Sigma h),$$
which yields
$$\nu^\prime(0)h=-\nabla_\Sigma h.$$
This in turn implies for the projection $P(\rho):=P_\Gamma$
$$P^\prime(\rho)h = -\nu^\prime(\rho)h\otimes \nu(\rho) -\nu (\rho)\otimes \nu^\prime(\rho)h,$$
hence
\begin{equation*}
P^\prime(0)h=\nabla_\Sigma h\otimes\nu_\Sigma + \nu_\Sigma \otimes \nabla_\Sigma h
=[\nabla_\Sigma\otimes \nu_\Sigma +\nu_\Sigma\otimes \nabla_\Sigma]h.
\end{equation*}
Applying these relations to $\nabla(\rho):=\nabla_\Gamma=P(\rho)M_0(\rho)\nabla_\Sigma$ yields
\begin{equation*}
(\nabla^\prime(0)h)\varphi
= (\nabla_\Sigma h|\nabla_\Sigma\varphi)\nu_\Sigma + hL_\Sigma\nabla_\Sigma\varphi
=[\nu_\Sigma\otimes \nabla_\Sigma h +h L_\Sigma ]\nabla_\Sigma \varphi,
\end{equation*}
and for a not necessarily tangent vector field $f$
$$(\nabla^\prime(0)h)f= \nu_\Sigma\otimes(\nabla_\Sigma h|\nabla_\Sigma)f+hL_\Sigma\nabla_\Sigma f.$$
For the divergence of the vector field $f$ this implies
$$[{\rm div}^\prime(0)h] f = (\nu_\Sigma|(\nabla_\Sigma h|\nabla_\Sigma)f)+h\, {\rm tr} [L_\Sigma\nabla_\Sigma f].$$
Finally,  the variation of the Laplace-Beltrami operator $\Delta(\rho):=\Delta_\Gamma$ becomes
$$(\Delta^\prime(0)h)\varphi= h\,{\rm tr}[L_\Sigma\nabla_\Sigma^2\varphi+\nabla_\Sigma(L_\Sigma\nabla_\Sigma\varphi)]+2(L_\Sigma\nabla_\Sigma h|\nabla_\Sigma\varphi)-\kappa(\nabla_\Sigma h|\nabla_\Sigma\varphi).$$
Note that in local coordinates we have
$${\rm tr}[L_\Sigma\nabla_\Sigma^2\varphi]= l^{ij}_\Sigma(\partial_i\partial_j\varphi -\Lambda_{ij}^k\partial_k\varphi),$$
hence with
$${\rm tr}[\nabla_\Sigma(L_\Sigma\nabla_\Sigma\varphi)]={\rm tr}[L_\Sigma\nabla_\Sigma^2\varphi]+({\rm div}_\Sigma L_\Sigma|\nabla_\Sigma\varphi),$$
we may write alternatively
$$(\Delta^\prime(0)h)\varphi= 2h\,{\rm tr}[L_\Sigma\nabla_\Sigma^2\varphi]+ (h\,{\rm div}_\Sigma L_\Sigma+[2L_\Sigma-\kappa_\Sigma]\nabla_\Sigma h|\nabla_\Sigma\varphi).$$
If $T=T_{ij}\tau^i\otimes\tau^j$ is a tensor we define
\begin{equation*}
{\rm div}_\Sigma T=(\tau^i|\partial_i(T_{ij}\tau^i))\tau^j+(\tau^i|T_{ij}\tau^i)\partial_i\tau^j.
\end{equation*}


\subsection{The Weingarten tensor and the mean curvature of $\Gamma$}
In invariant formulation we have with $P_\Gamma=:P(\rho)$
$$ L(\rho):=L_\Gamma= -\nabla_\Gamma \nu^\Gamma 
=- P(\rho) M_0(\rho) \nabla_\Sigma \{ \beta(\rho)(\nu_\Sigma -M_0(\rho)\nabla_\Sigma\rho)\}.$$
Thus for the variation of $L_\Gamma$ at $\rho=0$ we obtain with $P(0)=P_\Sigma$, $\beta(0)=1$, $M_0(0)=I$, 
and $P^\prime(0)=\nabla_\Sigma\otimes\nu_\Sigma \nu_\Sigma\otimes\nabla_\Sigma$,  
$\beta^\prime(0)=0$, $M_0^\prime(0)=L_\Sigma$,
$$ L^\prime(0)=  \nu_\Sigma\otimes L_\Sigma\nabla_\Sigma + L_\Sigma^2 +\nabla_\Sigma^2.$$
In particular, for $\kappa(\rho):=\kappa_\Gamma$ we have
$$\kappa(\rho) = - {\rm tr}[\nabla_\Gamma\nu^\Gamma] =
-{\rm tr}[P(\rho)M_0(\rho)\nabla_\Sigma\{\beta(\rho)(\nu_\Sigma-M_0(\rho)\nabla_\Sigma\rho)\}],$$
hence
\begin{equation}
\label{linearized-MC}
\kappa^\prime(0)= {\rm tr}\, L_\Sigma^2 +\Delta_\Sigma.
\end{equation}
Let us take another look at the mean curvature $\kappa_\Gamma$.
By the relations $\tau^r_\Gamma =P_\Gamma M_0(\rho) \tau^r_\Sigma$ and
$\nu^\Gamma=\beta(\rho)(\nu_\Sigma-a(\rho))$ we obtain
\begin{align*}
\kappa(\rho) &= - (\tau^j_\Gamma|\partial_j \nu^\Gamma)
=-(P_\Gamma M_0(\rho)\tau^j_\Sigma|(\partial_j\beta(\rho)/\beta(\rho))\nu^\Gamma
+\beta(\rho)(\partial_j \nu_\Sigma-\partial_j a(\rho)))\\
&= \beta(\rho)(P_\Gamma M_0(\rho)\tau^j_\Sigma|L_\Sigma \tau_j^\Sigma +\partial_j a(\rho))\\
&=\beta(\rho)(M_0(\rho)\tau^j_\Sigma|L_\Sigma\tau_j^\Sigma +\partial_j a(\rho))
- \beta(\rho)(\nu^\Gamma|M_0(\rho)\tau^j_\Sigma)(\nu^\Gamma|L_\Sigma \tau_j^\Sigma
+\partial_j a(\rho)).
\end{align*}
Since $(M_0(\rho)\tau^j_\Sigma|L_\Sigma\tau_j^\Sigma)={\rm tr}[M_0(\rho)L_\Sigma]$ as well as
$$(M_0(\rho)\tau^j_\Sigma|\partial_j a(\rho))= {\rm tr}[M_0(\rho)\nabla_\Sigma a(\rho)],$$ and
$(\nu^\Gamma|M_0\tau^j_\Sigma)= -\beta(\rho) [M_0(\rho)a(\rho)]^j$, we obtain
\begin{align*}
\kappa_\Gamma&= \beta(\rho)\big\{ {\rm tr}\big[M_0(\rho)(L_\Sigma+\nabla_\Sigma a(\rho))\big]\\
&+\beta^2(\rho)\big[M_0(\rho)a(\rho)\big]^j\big[(\nu_\Sigma|\partial_j a(\rho))- (a(\rho)|\partial_j a(\rho))
- (a(\rho)|L_\Sigma \tau_j^\Sigma)\big]\big\}\\
&= \beta(\rho)\big\{ {\rm tr}\big[M_0(\rho)(L_\Sigma+\nabla_\Sigma a(\rho))\big] 
-\beta^2(\rho)( M_0(\rho)a(\rho)|\nabla_\Sigma a(\rho)a(\rho))\big\},
\end{align*}
as $(\nu_\Sigma|a(\rho))=0$ implies
$$(\nu_\Sigma|\partial_j a(\rho))=- (\partial_j\nu_\Sigma|a(\rho))=
(L_\Sigma\tau_j^\Sigma|a(\rho)).$$
This yields the final form for the mean curvature of $\Gamma$.
\beq
\label{curvgamma}
\kappa(\rho) = \beta(\rho)\big\{ {\rm tr}\big[M_0(\rho)(L_\Sigma+\nabla_\Sigma a(\rho))\big] 
-\beta^2(\rho)( M_0(\rho)a(\rho)|[\nabla_\Sigma a(\rho)]a(\rho))\big\}.
\eeq
Recall that $a(\rho)=M_0(\rho)\nabla_\Sigma\rho$.

We can write the curvature of $\Gamma$ in local coordinates in the following form.
$$ \kappa(\rho) = c^{ij}(\rho,\partial\rho) \partial_i\partial_j\rho
+g(\rho,\partial\rho),$$ with
$$c^{ij} = \beta(\rho) [M_0^2(\rho)]^{ij}-\beta(\rho)^3[M_0^2\nabla_\Sigma\rho]^i
[M_0^2\nabla_\Sigma\rho]^j.$$
A simple computation yields for the symbol $c(\rho,\xi)=c^{ij}\xi_i\xi_j$
of the principal part of this operator
$$c(\rho,\xi)= \beta(\rho)\{ |M_0(\rho)\xi|^2- \beta^2(\rho)(a(\rho)|M_0(\rho)\xi)^2\}\geq
\beta^3(\rho)|M_0(\rho)\xi|^2\geq \eta |\xi|^2,$$
for $\xi =\xi_k\tau^k_\Sigma\in T_p\Sigma$,
as long as $|\nabla_\Sigma\rho|_{\infty}<\infty$ and $|\rho|_\infty<\rho_0.$ Therefore
the curvature $\kappa(\rho)$ is a quasi-linear strongly elliptic differential operator
on $\Sigma$, acting on the parametrization $\rho$ of $\Gamma$ over $\Sigma$,
see also \cite{ES97,ES98} for a different derivation.

\medskip
\goodbreak

\section{Approximation of hypersurfaces} 
\subsection{The tubular neighborhood of a hypersurface}
Let $\Sigma$ be a compact connected $C^2$-hypersurface bounding
a domain $\Omega\subset\RR^n$, and
let $\nu_\Sigma$ be the outer unit normal field on $\Sigma$
with respect to $\Omega$.
\\ 
The conditions imply that $\Sigma$  satisfies a {\em uniform interior and exterior ball condition}, 
i.e.\ there is a number $a>0$ such that for each 
point $p\in\Sigma$ there are balls $B_a(x_i)\subset \Omega_i$ such that $\Sigma\cap \bar{B}_a(x_i)=\{p\}$. 
As in \cite[Section 14.6]{GT01} we conclude that the mapping
\begin{equation}
\label{X}
X :\Sigma\times (-a,a)\to \RR^{n},\qquad X(p,r):=p+r\nu_\Sigma(p) 
\end{equation}
is a $C^1$-diffeomorphism onto its image $ U_a:=\text{im}(X)$.
It will be convenient to decompose the inverse
of $X$ into $X^{-1}=(\Pi_\Sigma,d_\Sigma)$ such that
$$\Pi_\Sigma\in C^{1}(U_a,\Sigma),
\qquad d_\Sigma\in C^{1}(U_a,(-a,a)).$$
$\Pi_\Sigma(x)$ is the nearest point on $\Sigma$ to $x$,
$d_\Sigma(x)$ is the signed distance from $x$ to $\Sigma$, and
 $U_a$ consists of  the set of those points in $\RR^n$ which have
distance less than $a$ to $\Sigma$, and $|d_\Sigma(x)|=dist(x,\Sigma)$,
$d_\Sigma(x)<0$ if and only if $x\in\Omega$.

\medskip

{(i)} From the uniform interior and exterior ball condition follows that
the number $1/a$ bounds the principal curvatures of $\Sigma$, i.e.,
\begin{equation}
\label{a-curvatures}
\max\{\kappa_i(p): p\in\Sigma,\ i=1,\cdots ,n-1\}\le 1/a.
\end{equation}

{(ii)} We remark here that the regularity assertion
$X^{-1}\in C^1(U_a,\Sigma\times (-a,a))$ is an easy consequence of
the inverse function theorem.
To see this, we fix a point $(p_0, r_0)\in \Sigma\times (-a,a)$ and a chart $\phi$ for $p_0$.
Then the function $f(\theta,r)=X(\phi(\theta),r)$ has derivative
\begin{equation*}
Df(0,r_0)=[[I-r_0L_\Sigma(p_0)]\phi^\prime(0),\nu_\Sigma(p_0)].
\end{equation*}
It follows from \eqref{a-curvatures} that
$[I-r_0L_\Sigma(p_0)]\in\cB(T_{p_0}\Sigma)$ is invertible,
and consequently, $Df(0,r_0)\in\cB(\R^n)$ is invertible as well. 
The inverse function theorem implies that $X$ is locally invertible with inverse
of class $C^1$. In particular, $\Pi_\Sigma$ and $d_\Sigma$ are $C^1$.

\medskip

{(iii)}  A remarkable fact is that the signed distance $d_\Sigma$ is even of class $C^2$. To see this, we use the identities
$$x-\Pi_\Sigma(x)=d_\Sigma(x)\nu_\Sigma(\Pi_\Sigma(x)),\quad d_\Sigma(x) = (x-\Pi_\Sigma(x)|\nu_\Sigma(\Pi_\Sigma(x)).$$
Differentiating w.r.t.\ $x_k$ this yields
\begin{align*}
\partial_{x_k} d_\Sigma(x)&= (e_k-\partial_{x_k} \Pi_\Sigma(x)|\nu_\Sigma(\Pi_\Sigma(x)))+
(x-\Pi_\Sigma(x)|\partial_{x_k} (\nu_\Sigma\circ \Pi_\Sigma)(x))\\
 &= \nu_k(\Pi_\Sigma(x)) +d_\Sigma(x)(\nu_\Sigma(\Pi_\Sigma(x))|\partial_{x_k}(\nu_\Sigma\circ\Pi_\Sigma(x)) )
= \nu_k(\Pi_\Sigma(x)),
\end{align*}
since $\partial_{x_k} \Pi_\Sigma(x)$ belongs to the tangent space $T_{\Pi_\Sigma(x)}\Sigma$, as does
$\partial_{x_k}(\nu_\Sigma\circ \Pi_\Sigma(x))$, since $|\nu_\Sigma\circ\Pi_\Sigma(x)|\equiv1$.
Thus we have the formula
\begin{equation}\label{dd}
\nabla_x d_\Sigma(x)= \nu_\Sigma(\Pi_\Sigma(x)),\quad x\in U_a.
\end{equation}
This shows, in particular, that $d_\Sigma$ is of class $C^2$.

\medskip

{(iv)}\
It is useful to also have a representation of $\nabla_x \Pi_\Sigma(x)$. With 
$$ I-\Pi_\Sigma^\prime(x) = d_\Sigma^\prime(x) \nu_\Sigma(\Pi_\Sigma(x))+d_\Sigma(x) \nu_\Sigma^\prime(\Pi_\Sigma(x))\Pi^\prime_\Sigma(x),$$
and \eqref{dd}, we obtain
\begin{equation}\label{dPi}\nabla_x\Pi_\Sigma(x)= P_\Sigma (\Pi_\Sigma(x)) M_0(d_\Sigma(x))(\Pi_\Sigma(x)),
\end{equation}
where $M_0(r)(p):=(I-rL_\Sigma(p))^{-1}$.
This shows, in particular, that
that $\nabla_x\Pi_\Sigma(p)=P_\Sigma(p)$ is the orthogonal projection onto the tangent space $T_p\Sigma$.


\subsection{The level function}
Let $\Sigma$ be a compact connected hypersurface of class $C^2$ bounding the domain $\Omega$ in $\RR^n$. According
to the previous section, $\Sigma$ admits a tubular neighborhood $U_a$ of width $a>0$. We may assume w.l.o.g.\ $a\leq 1$. The signed distance function $d_\Sigma(x)$
in this tubular neighborhood is of class $C^2$ as well, and since
$$\nabla_x d_\Sigma(x) = \nu_\Sigma(\Pi_\Sigma(x)),\quad  x\in U_a,$$
we can view $\nabla_x d_\Sigma(x)$ as a $C^1$-extension of the normal field $\nu_\Sigma(x)$ from $\Sigma$ to the
tubular neighborhood $U_a$ of $\Sigma$.
Computing the second derivatives $\nabla_x^2 d_\Sigma$ we obtain
\begin{equation*}
\begin{split}
\nabla^2_x d_\Sigma(x)=\nabla_x \nu_\Sigma (\Pi_\Sigma(x))
&=-L_\Sigma(\Pi_\Sigma(x))P_\Sigma (\Pi_\Sigma(x))(I-d_\Sigma(x)L_\Sigma(\Pi_\Sigma(x)))^{-1}\\
&= -L_\Sigma(\Pi_\Sigma(x))
(I-d_\Sigma(x)L_\Sigma(\Pi_\Sigma(x)))^{-1},
\end{split}
\end{equation*}
for $x\in U_a$, as $L_\Sigma(p)=L_\Sigma (p) P_\Sigma(p)$.
Taking traces then yields
\begin{equation}
\label{D2-distance}
\Delta d_\Sigma(x)=-\sum_{i=1}^{n-1}\frac{\kappa_i(\Pi_\Sigma(x))}{1-d_\Sigma(x)\kappa_i(\Pi_\Sigma(x))},
\quad x\in U_a.
\end{equation}
In particular, this implies 
\begin{equation}
\label{curvature-level}
\kappa_\Sigma (p)=-\Delta_x d_\Sigma(p),\quad p\in\Sigma.
\end{equation} 

Therefore the norm of $\nabla_x^2 d_\Sigma$ is equivalent to the maximum of the moduli of
the curvatures  of $\Sigma$ at a fixed point.
Hence we find a constant $c$, depending only on $n$, such that
$$ c|\nabla_x^2 d_\Sigma|_\infty\leq\max\{|\kappa_i(p)|:\, i=1,\ldots,n-1,\, p\in\Sigma\}
\leq c^{-1}|\nabla_x^2 d_\Sigma|_\infty.$$
It has now become clear that the Lipschitz constant for the normal $\nu_\Sigma(p)$,
which is given by $|\nabla_x^2 d_\Sigma|_{\infty}$,
is equivalent to the maximum of the moduli of the principal curvatures of $\Sigma$.

Next we extend $d_\Sigma$ as a function $\varphi$ to all of $\RR^n$. For this purpose we choose a $C^\infty$-function $\chi(s)$ such that
$\chi(s)=1$ for $|s|\leq 1$, $ \chi(s)= 0$ for $|s|\geq 2$, $0\leq \chi(s)\leq1$. Then we set
\begin{equation}
\label{level}
\varphi(x)= d_\Sigma(x)\chi(3d_\Sigma(x)/a)+ ({\rm sign}\, d_\Sigma(x))(1-\chi(3d_\Sigma(x)/a)),\quad x\in U_a,
\end{equation}
and $\varphi=1$ in the exterior component of $\RR^n\setminus U_a$,  $\varphi=-1$ in its interior component.
This function $\varphi$ is then of class $C^2$, 
$\varphi(x)=d_\Sigma(x)$ for $x\in U_{a/3}$, and
$$
\text{$\varphi(x)=0$ if and only if $x\in\Sigma.$}
$$
Thus $\Sigma$ is the
level set $\Sigma=\varphi^{-1}(0)$ of $\varphi$ at level $0$,
$\varphi$ is called a {\em canonical level function} for $\Sigma$. It is a special level function for $\Sigma$, as $\nabla_x\varphi(x)=\nu_\Sigma(x)$
for $x\in\Sigma$.


\subsection{Existence of parameterizations}
Recall the Haussdorff metric on the set $\mathcal K$ of compact subsets of $\RR^n$ defined by
$$ d_H(K_1,K_2)=\max\{ \sup_{x\in K_1} d(x,K_2), \sup_{x\in K_2} d(x,K_1)\}.$$
Suppose $\Sigma$ is a compact (and without loss of generality) connected hypersurface of class $C^2$ in $\RR^n$. As before, let $U_a$ be its tubular neighborhood, $\Pi_\Sigma:U_a\rightarrow\Sigma$ the projection and
$d_\Sigma:U_a\rightarrow\RR$ the signed distance. We want to parameterize hypersurfaces $\Gamma$ which are close to $\Sigma$
as 
$$q=p+\rho(p)\nu_\Sigma(p),$$ 
where
$\rho:\Sigma\rightarrow\RR$ is then called the {\it normal parametrization} of $\Gamma$ over $\Sigma$. For this to make sense, $\Gamma$ must belong to the tubular neighborhood
$U_a$ of $\Sigma$. Therefore, a natural requirement would be $d_H(\Gamma,\Sigma)<a$. We then say that $\Gamma$ and $\Sigma$ are {\em $C^0$-close} (of order $\varepsilon$) if  $d_H(\Gamma,\Sigma)<\varepsilon$.

However, this condition is not enough to allow for existence of the parametrization, since it is not clear that
the map $\Pi_\Sigma$ is injective on $\Gamma$:
small Haussdorff distance does not prevent $\Gamma$ from folding within the tubular neighborhood. We need a stronger assumption to prevent this. If $\Gamma$ is
a hypersurface of class $C^1$ we may introduce the so-called {\em normal bundle} $\cN\Gamma$ defined by
$$\cN\Gamma:=\{(p,\nu_\Gamma(p)):\, p\in\Gamma\}\subset\RR^{2n}.$$
Suppose $\Gamma$ is a compact, connected $C^1$-hypersurface in $\R^n$. 
We say that $\Gamma$ and $\Sigma$ are {\em $C^1$-close} (of order $\varepsilon$) if $d_H(\cN\Gamma,\cN\Sigma)<\varepsilon$.
We are going to show that
$C^1$-hypersurfaces $\Gamma$ which are $C^1$-close to $\Sigma$ can in fact be parametrized over $\Sigma$.

For this purpose observe that, in case $\Gamma$ and $\Sigma$ are $C^1$-close of order $\varepsilon$, whenever $q\in\Gamma$, then there is $p\in\Sigma$ such that
$|q-p|+|\nu_\Gamma(q)-\nu_\Sigma(p)|<\varepsilon$.
Hence $|q-\Pi_\Sigma q|<\varepsilon$, with $\Pi_\Sigma q:=\Pi_\Sigma(q)$, and
$$|\nu_\Gamma(q)-\nu_\Sigma(\Pi_\Sigma q)|\leq |\nu_\Gamma(q)-\nu_\Sigma(p)|+|\nu_\Sigma(\Pi_\Sigma q)-\nu_\Sigma(p)|\leq \varepsilon+L|\Pi_\Sigma q-p|,$$
which yields with $|\Pi_\Sigma q-p|\leq|\Pi_\Sigma q-q|+|p-q|<2\varepsilon$,
$$|q-\Pi_\Sigma q|+|\nu_\Gamma(q)-\nu_\Sigma(\Pi_\Sigma q)|\leq 2(1+L)\varepsilon,$$
where $L$ denotes the Lipschitz constant of the normal of $\Sigma$. In particular, the tangent space $T_{q}\Gamma$ is transversal
to $\nu_\Sigma(\Pi_\Sigma q)$, for each $q\in\Gamma$, that is, 
$$
T_q\Gamma \oplus {\rm span\,}\{\nu_\Sigma(\Pi_\Sigma q)\}=\R^n,\quad q\in\Gamma.
$$
Now fix a point  $q_0\in\Gamma$ and set $p_0=\Pi_\Sigma q_0$. Since the tangent space $T_{q_0}\Gamma$ is transversal
to $\nu_\Sigma(p_0)$, we see that
$\Pi^\prime_\Sigma(q_0):T_{q_0}\Gamma\rightarrow T_{p_0}\Sigma$ is invertible. The inverse function theorem yields an open neighborhood
$V(p_0)\subset\Sigma$ and a
$C^1$-map $g:V(p_0)\rightarrow\Gamma$ such that $g(p_0)=q_0$, $g(V(p_0))\subset\Gamma$, and $\Pi_\Sigma g(p)=p$ in $V(p_0)$. Therefore we obtain
\begin{equation*}
q=g(p)=\Pi_\Sigma g(p)+d_\Sigma(g(p))\nu_\Sigma(\Pi_\Sigma g(p))=p+\rho(p)\nu_\Sigma(p), \quad \rho(p):=d_\Sigma(g(p)).
\end{equation*}
Thus we have a local parametrization of $\Gamma$ over $\Sigma$. 
We may extend $g$ to a
maximal domain $V\subset\Sigma$, e.g.\ by means of Zorn's lemma. 
Clearly $V$ is open in $\Sigma$ and we claim that $V=\Sigma$. 
If not, then the boundary of $V$ in $\Sigma$ is nonempty and
hence we find a sequence $p_n\in V$ such that $p_n\rightarrow p_\infty\in\partial V$.
Since $\rho_n=\rho(p_n)$ is bounded, we may assume w.l.o.g. that 
$\rho_n\rightarrow\rho_\infty$. But then $q_\infty=p_\infty+\rho_\infty\nu_\Sigma(p_\infty)$
belongs to $\Sigma$ as $\Sigma$ is closed. 
Now we may apply the inverse function theorem again to see that $V$ cannot be maximal.
Since the map
$\Phi(p)=p+\rho(p)\nu_\Sigma(p)$ is a local $C^1$-diffeomorphism, it is also open. 
Hence $\Phi(\Sigma)\subset\Gamma$ is open and compact, i.e.\ $\Phi(\Sigma)=\Gamma$
by connectedness of $\Gamma$. The map $\Phi$ is therefore a $C^1$-diffeomorphism from $\Sigma$ to $\Gamma$.
In case $\Sigma$ is of class $C^3$ the proof above immediately implies that
$\Phi\in {\rm Diff}^2(\Sigma,\Gamma)$.

Observe that because of $x=\Pi_\Sigma x+d_\Sigma(x)\nu_\Sigma(\Pi_\Sigma x)$ in $U_a$ we have $x\in\Gamma$ 
if and only if $d_\Sigma(x)=\rho(\Pi_\Sigma x)$. This property can be used to construct
a $C^1$-function $\psi$ on $\RR^n$ such that $\Gamma=\psi^{-1}(0)$, i.e.\ a level function for $\Gamma$.
For example we may take
$$\psi(x)=\varphi(x)-\rho(\Pi_\Sigma x)\chi(3d_\Sigma(x)/a),\quad x\in\RR^n,$$
provided $\varepsilon<a/3$, where $\varphi$ and $\chi$ are as in subsection 2.


\subsection{Approximation of hypersurfaces}
Suppose as before that $\Sigma$ is a compact connected hypersurface of class $C^2$ enclosing a domain $\Omega$ in 
$\R^n$.
We may use the level function $\varphi:\RR^n\rightarrow\RR$ 
introduced in \eqref{level} 
to construct a real analytic hypersurface $\Sigma_\eps$
such that $\Sigma$ appears as a $C^2$-graph over $\Sigma_\eps$.
In fact, we  show that there is  $\eps_0\in (0,a/3)$ 
such that for every $\eps\in (0,\eps_0)$ there is an analytic manifold $\Sigma_\eps$ and 
a function $\rho_\eps\in C^2(\Sigma_\eps)$ with the property that
$$\Sigma=\{p+\rho_\eps(p)\nu_{\Sigma_\eps}(p): p\in \Sigma_\eps\}$$
and
$$|\rho_\eps|_\infty+|\nabla_{\Sigma_\eps}\rho_\eps|_\infty
+|\nabla^2_{\Sigma_\eps}\rho_\eps|_\infty\leq\eps.  $$
For this purpose, choose $R>0$ such that $\varphi(x)=1$ for $|x|>R/2$. Then define
$$\psi_k(x)=c_k\Big(1-\frac{|x|^2}{R^2}\Big)^k_+,\quad x\in\RR^n,$$
where $c_k>0$ is chosen such that $\int_{\RR^n} \psi_k(x)dx=1$; note that $c_k\sim k^{n/2}$ as $k\to\infty$.
Then as $k\rightarrow\infty$, we have $\psi_k(x)\rightarrow0$, uniformly
for $|x|\geq\eta>0$, hence $\psi_k*f\rightarrow f$ in $BU\!C^m(\RR^n)$, whenever $f\in BU\!C^m(\RR^n)$. We define
$\varphi_k=1+\psi_k*(\varphi-1)$; then $\varphi_k\rightarrow\varphi$ in $BU\!C^2(\RR^n)$. Moreover,
$$\psi_k*(\varphi-1)(x)=\int_{\RR^n} (\varphi(y)-1)\psi_k(x-y)dy=\int_{B_{R/2}(0)} (\varphi(y)-1)\psi_k(x-y)dy.$$
For $|x|, |y|<R/2$ follows $|x-y|<R$, and hence $\psi_k(x-y)=c_k(1-|x-y|^2/R^2)^k$ is  polynomial in $x,y$. But then
$\varphi_k(x)$ is a polynomial 
for such values of $x$, in particular $\varphi_k$ is real analytic in $U_a$. Choosing $k$ large enough,
we have $|\varphi-\varphi_k|_{BU\!C^2(\RR^n)}<\eps$.

Now suppose $\varphi_k(x)=0$; then $|\varphi(x)|<\eps$, hence
$x\in U_a$ and then $|d_\Sigma(x)|<\eps$. This shows that the set $\Sigma_k=\varphi_k^{-1}(0)$ is in the $\eps$-tubular neighborhood
around $\Sigma$. Moreover, $|\nabla\varphi_k-\nabla\varphi|_\infty<\eps$ yields $\nabla\varphi_k(x)\neq0$ in $U_a$, and therefore
$\Sigma_k$ is a manifold, which is real analytic.

Next we show that $\Sigma$ and $\Sigma_k$ are $C^1$-diffeomorphic. For this purpose, fix a point $q_0\in\Sigma_k$.
Then $q_0=p_0+\ro_0\nu_\Sigma(p_0)$, where $p_0=\Pi_\Sigma q_0\in\Sigma$ and $\ro_0=d_\Sigma(q_0)$. Consider the equation
$g(p,\ro):=\varphi_k(p+\ro\nu_\Sigma(p))=0$ near $(p_0,\ro_0)$. 
Since 
$$
\partial_\ro g(p,\ro)=
(\nabla_x\varphi_k(p+\ro\nu_\Sigma(p))|\nu_\Sigma(p))$$ 
we have 
\begin{align*}
\partial_\ro g(p_0,\ro_0)&=(\nabla\varphi_k(q_0)|\nabla\varphi(p_0))\\
&\geq 1-|\nabla\varphi_k(q_0)-\nabla\varphi(q_0)|-|\nabla\varphi(q_0)-\nabla\varphi(p_0)|\\
&\geq 1-|\varphi_k-\varphi|_{BC^1(\RR^n)}- a|\nabla ^2\varphi|_{BC(\RR^n)}\geq 1-\eps-aL>0.
\end{align*}
Therefore, we may apply the implicit function theorem to obtain an open neighborhood $V(p_0)\subset\Sigma$ 
and a $C^1$-function $\ro_k:V(p_0)\rightarrow \RR$
such that $\ro_k(p_0)=\ro_0$ and $p+\ro_k(p)\nu_\Sigma(p)\in\Sigma_k$ for all $p\in V(p_0)$. 
We can now proceed as in subsection 3 to extend $\ro_k(\cdot)$ to a maximal domain $V\subset\Sigma$,
which coincides with $\Sigma$ by compactness and connectedness of $\Sigma$. 

Thus we have a well-defined $C^1$-map $f_k:\Sigma\rightarrow\Sigma_k$, $f_k(p)=p+\ro_k(p)\nu_\Sigma (p)$,
which is injective and a diffeomorphism from $\Sigma$ to its range. 
We claim that $f_k$ is also surjective. If not, there is some point
$q\in\Sigma_k $, $q\not\in f_k(\Sigma)$. Set $p=\Pi_\Sigma q$. 
Then $q=p+d_\Sigma(p)\nu_\Sigma(p)$ with $d_\Sigma(p)\neq \ro_k(p)$.
Thus, there there are at least two numbers $\ro_1,\ro_2\in (-a,a)$
with $p+\ro_i\nu_\Sigma(p)\in\Sigma_k$. This implies with $\nu_\Sigma=\nu_\Sigma(p)$
$$0=\varphi_k(p+\ro_2\nu_\Sigma)-\varphi_k(p+\ro_1\nu_\Sigma)
=(\ro_2-\ro_1)\int_0^1(\nabla\varphi_k(p+(\ro_1+ t(\ro_2-\ro_1))\nu_\Sigma)|\nu_\Sigma)\,dt,$$
which yields $\ro_2-\ro_1=0$ since
$$\int_0^1(\nabla\varphi_k(p+(\ro_1+ t(\ro_2-\ro_1))\nu_\Sigma)|\nu_\Sigma)\,dt\geq 1-\eps-aL>0,$$
as above. Therefore the map $f$ is also surjective, and hence $f_k\in {\rm Diff}^1(\Sigma,\Sigma_k)$.
This implies in particular that $\Sigma_k=f_k(\Sigma)$ is connected.
For later use we note that 
$$|\ro_k|_\infty+|\nabla_\Sigma \ro_k|_\infty \to 0\;\;\text{as}\;\; k\to \infty, $$
as can be inferred from the relationship $\varphi_k(p+r_k(p)\nu_\Sigma(p))=0$ for $p\in\Sigma$.

\medskip

Next we show that the mapping
\begin{equation*}
\label{X-k-diff}
X_k:\Sigma_k\times (-a/2,a/2)\to U(\Sigma_k,a/2),\quad X_k(q,s):=q+s\nu_k(q)
\end{equation*} 
is a $C^1$-diffeomorphism for $k\ge k_0$, with $k_0\in \BN$ sufficiently large.
In order to see this, we use the diffeomorphism $f_k$ constructed above to rewrite $X_k$ as
\begin{equation*}
\begin{split}
X_k(q,s)&=X_k(f_k(p),s)\\
&=p + s\,\nu_\Sigma(p) +\ro_k(p) \nu_\Sigma(p)+ s[\nu_k(p+\ro_k(p)\nu_\Sigma(p))-\nu_\Sigma(p)]\\
&=:X(p,s)+G_k(p,s)=:H_k(p,s).
\end{split}
\end{equation*}
Clearly $H_k\in C^1(\Sigma\times (-a/2,a/2),\R^n)$ and 
$X\in {\rm Diff}^1(\Sigma\times (-a,a), U(\Sigma,a))$. 
It is not difficult to see that
\begin{equation*}
|G_k(p,s)|+|DG_k(p,s)|\to 0\;\; \text{as} \;\; k\to \infty, \;\;
\text{uniformly in $(p,s)\in \Sigma\times [-a/2,a/2]$.}
\end{equation*}
Consequently, $DH_k(p,s):T_p(\Sigma)\times (-a/2,a/2)\to \R^n$ is invertible for $k\ge k_0$,
and by the inverse function theorem, $H_k$ is a local $C^1$-diffeomorphism.
We claim that $H_k$ is injective for all  $k$ sufficiently large.
For this purpose, note that due to compactness of $\Sigma\times [-a/2,a/2]$ and injectivity of $X$ there exists 
a constant $c>0$ such that
\begin{equation*}
\label{X-estimate}
|X(p,s)-X(\bar p,\bar s)|\ge c\big(|p-\bar p|+|s-\bar s|\big),
\quad (p,s),\; (\bar p,\bar s)\in \Sigma\times [-a/2,a/2].
\end{equation*}
The properties of $G_k$ and compactness of $\Sigma\times [-a/2,a/2]$ imply, in turn, that
the estimate above remains true for $X$ replaced by $H_k$, 
and  $c$ replaced by $c/2$, provided $k\ge k_0$ with $k_0$ sufficiently large.
Hence $H_k$ is a $C^1$-diffeomorphism onto its image for $k$ sufficiently large, as claimed.
This shows that  $\Sigma_k$ has a uniform tubular neighborhood of width $a/2$ for any $k\ge k_0$,
and it follows that $\Sigma\subset U(\Gamma_k,a/2)$.
$\Sigma$ and $\Sigma_k$ are compact, connected, $C^1$ hypersurfaces,
and may now apply the results of subsection 3, showing that $\Sigma$ can be parameterized over $\Sigma_k$ by means of
$$ [p\mapsto p+\rho_k(p)\nu_{k}(p)]\text{ where }\rho_k\in C^2(\Sigma_k,\R).$$
Finally, it is not difficult to show that the relation $\varphi(p+\rho_k(p)\nu_k(p))=0$ for $p\in\Sigma_k$ implies
$
|\rho_k|_{\infty}+|\nabla_{\Sigma_k}\rho_k|_\infty
+|\nabla^2_{\Sigma_k}\rho_k|_\infty\leq \eps 
$
for $k$ sufficiently large. 
\goodbreak
\section{Compact embedded hypersurfaces in $\R^n$}
\subsection{The manifold of compact connected hypersurfaces of class $C^{2}$}
Consider the set $\mathcal{M}$ of all compact connected $C^{2}$-hypersurfaces $\Sigma$
in $\RR^n$.
Let $\cN\Sigma$
denote their associated normal bundles.
The second normal bundle of $\Sigma$ is defined by
$$\cN^2\Sigma=\{(p,\nu_\Sigma(p),\nabla_\Sigma\nu_\Sigma(p)):\, p\in\Sigma\}.$$
We introduce a metric $d_{\mathcal{M}}$ on $\mathcal{M}$ by means of
$d_{\mathcal{M}}(\Sigma_1,\Sigma_2)=d_H(\cN^2\Sigma_1,\cN^2\Sigma_2)$.
This way $\mathcal{M}$ becomes a metric space. We want to show  that $\mathcal{M}$ is
a Banach manifold.

Fix a hypersurface $\Sigma\in\mathcal{M}$ of class $C^3$. Then we define a
chart over the Banach space $X_\Sigma:=C^2(\Sigma,\RR)$ as follows.
$\Sigma$ has a tubular neighborhood $U_a$ of width $a$. Therefore we take as the chart set,
say $B^\Sigma_{a/3}(0)\subset X_\Sigma$, and for
a given function $\rho\in B^{\Sigma}_{a/3}(0)$, we define the hypersurface
$\Gamma^\Sigma_\rho$ by means of the map
$$\Phi^\Sigma(\rho)(p)=p+\rho(p)\nu_\Sigma(p),\quad p\in\Sigma.$$
According to Section 4, this yields a hypersurface $\Gamma_\rho^\Sigma$ of class $C^2$,
diffeomorphic to $\Sigma$. Moreover, with some constant $C^\Sigma_a$,
we have
$$d_{\mathcal{M}}(\Gamma^\Sigma_\rho,\Sigma)\leq C^\Sigma_a |\rho|_{BU\!C^2(\Sigma)},$$
which shows that the map $\Phi^\Sigma(\rho):B^\Sigma_{a/2}(0)\rightarrow \mathcal{M}$
is continuous. Conversely, given  $\Gamma\in\mathcal{M}$
which is $C^2$-close to $\Sigma$, the results in subsection 4.3 show that
$\Gamma$ can be parameterized by a function $\rho\in C^2(\Sigma,\RR)$, such that
$|\rho|_{BU\!C^2(\Sigma)}<a/3$.

We compute the tangent space $T_\Sigma\cM$ at some fixed $\Sigma\in\mathcal{M}$. For this purpose
we take a differentiable curve
$\Gamma:(-\delta_0,\delta_0)\rightarrow \mathcal{M}$ with $\Gamma(0)=\Sigma$.
Then according to subsection 4, there is $\delta\in(0,\delta_0)$ such that
for each $t\in(-\delta,\delta)$ we find a parametrization $\rho(t)\in C^2(\Sigma,\RR)$
of $\Gamma(t)$. Then in these coordinates we have
$$V:=\frac{d}{dt}\Gamma(t)\Big|_{t=0}=\frac{d}{dt}\rho(0)\in C^2(\Sigma,\RR)=X_\Sigma.$$
In  other words, the tangent space $T_\Sigma\mathcal{M}$ consists of
all normal velocity fields $V$ on $\Sigma$ which are of class $C^2$.
Moreover, if $\psi(t,x)$ is the level  function for $\Gamma(t)$ from subsection 4.2, then
$$0=\psi\big(t,\phi(\theta)+\rho(t,\phi(\theta))\nu_\Sigma(\phi(\theta))\big),$$
hence at $t=0$
$$0=\partial_t \psi + (\nabla_x\psi|\partial_t\rho\nu_\Sigma)=
\partial_t\psi + \partial_t\rho(\nu_\Sigma|\nu_\Sigma) = \partial_t\psi +V_\Gamma,$$
i.e.\ we have $\partial_t\psi= -V_\Gamma$, for the normal velocity $V_\Gamma$ of
$\Gamma(t)$.

There is one shortcoming with this approach, namely the need to require that
$\Sigma\in C^3$. This is due to the fact that we are loosing one derivative when
forming the normal $\nu_\Sigma$. However, since we may approximate a given hypersurface
of class $C^2$ by a real analytic one in the second normal bundle, this defect can be
avoided by only parameterizing over real analytic hypersurfaces which is sufficient.


\subsection{Compact hypersurfaces with uniform ball condition.}
Let $\Omega\subset\RR^n$ be a bounded domain and consider a closed connected $C^2$-hypersurface $\Gamma\subset\Omega$. This hypersurface separates $\Omega$ into two disjoint open connected sets $\Omega_1$ and $\Omega_2$, the interior and the exterior of $\Gamma$ w.r.t.\ $\Omega$. By means of the level function $\varphi_\Gamma$ of $\Gamma$ 
we have 
$\Omega_1=\varphi_\Gamma^{-1}(-\infty,0)$ and $\Omega_2=\Omega\setminus\bar{\Omega}_1$. Then $\partial\Omega_1=\Gamma$ and $\partial\Omega_2 =\partial\Omega\cup \Gamma$.

The hypersurface $\Gamma$  satisfies the {\em ball condition}, i.e.\ there is a radius $r>0$ such that for each point $p\in\Gamma$ there are balls $B_r(x_i)\subset \Omega_i$ such that $\Gamma\cap \bar{B}_r(x_i)=\{p\}$. The set of  hypersurfaces of class $C^2$ contained in $\Omega$ satisfying
the ball condition with radius $r>0$ will be denoted by $\cM^2(\Omega,r)$. Note that hypersurfaces in this class have uniformly bounded principal curvatures.

The elements of $\cM^2(\Omega,r)$ have a tubular neighborhood of width $a$  larger than $r/2$. Therefore the construction of the level function $\varphi_\Gamma$ of $\Gamma $ from subsection 4.2 can be carried out with the same $a$ and the same cut-off function $\chi$ for each $\Gamma\in\cM^2(\Omega,r)$.
More precisely, we have
$$\varphi_\Gamma(x) = g(d_\Gamma(x)),\quad x\in\Omega,$$
with
$$g(s) =s\chi(3s/a) +{\rm sgn}(s)(1-\chi(3s/a)),\quad s\in\RR;$$

\medskip

\noindent
note that $g$ is strictly increasing and equals $\pm 1$ for $\pm d_\Gamma(x)>2a/3$.
This induces an injective map $\Phi:\cM^2(\Omega,r)\to C^2(\bar{\Omega})$ which assigns to $\Gamma$ the level function $\varphi_\Gamma$. $\Phi$ is in fact an isomorphism of  $\cM^2(\Omega,r)$ onto $\Phi(\cM^2(\Omega,r))\subset C^2(\bar{\Omega})$.

This can be seen as follows; let $\eps>0$ be small enough. If $|\varphi_{\Gamma_1}-\varphi_{\Gamma_2}|_{2,\infty}\leq\eps$ then $d_{\Gamma_1}(x)\leq \eps$ on $\Gamma_2$ and
$d_{\Gamma_2}(x)\leq \eps$ on $\Gamma_1$, which implies $d_H(\Gamma_1,\Gamma_2)\leq \eps$.
Moreover, we also have $|\nabla_x \varphi_{\Gamma_1}(x)-\nu_{\Gamma_2}(x)|\leq \eps$ on $\Gamma_2$ and $|\nabla_x \varphi_{\Gamma_2}(x)-\nu_{\Gamma_1}(x)|\leq \eps$ on $\Gamma_1$ which yields $d_H(\cN\Sigma_1,\cN\Sigma_2)\leq C_0\eps$.
Then the hypersurfaces $\Gamma_j$ can both be parameterized over a $C^3$-hypersurface $\Sigma$, and therefore $d_H(\cN^2\Gamma_1,\cN^2\Gamma_2)\leq\eps$ if and only if
$$|\rho_1-\rho_2|_\infty +|\nabla_\Sigma(\rho_1-\rho_2)|_\infty+|\nabla_\Sigma^2(\rho_1-\rho_2)|_\infty<C_1\eps.$$
This in turn is equivalent to $|\varphi_{\Gamma_1}-\varphi_{\Gamma_2}|_{2,\infty}\leq C_2\eps.$

Let $s-(n-1)/p>2$; for $\Gamma\in\cM^2(\Omega,r)$, we define $\Gamma\in W^s_p(\Omega,r)$ if $\varphi_\Gamma\in W^s_p(\Omega)$. In this case the local charts for $\Gamma$ can be chosen of class $W^s_p$ as well. A subset $A\subset W^s_p(\Omega,r)$ is said to be (relatively) compact, if $\Phi(A)\subset W^s_p(\Omega)$
is (relatively) compact. In particular, it follows from Rellich's theorem that $W^s_p(\Omega,r)$ is a compact subset of $W^\sigma_q(\Omega,r)$, whenever  $s-n/p>\sigma-n/q$, and $s>\sigma$.


\medskip
Received xxxx 20xx; revised xxxx 20xx.
\medskip


\begin{thebibliography}{99}

\bibitem{DoC92} 
\newblock M.P.~Do Carmo, 
\newblock{``Riemannian Geometry,"}
\newblock Mathematics: Theory \& Applications, Birkh\"auser, Basel, 1992.

\bibitem{ES97} 
\newblock J. Escher and G. Simonett,
\newblock \emph{Classical solutions for Hele-Shaw models with surface tension},
\newblock {Adv. Differential Equations} \textbf{2} (1997), 619--642.
 
\bibitem{ES98} 
\newblock J. Escher and G. Simonett,
\newblock \emph{A center manifold analysis for the Mullins-Sekerka model},
\newblock {J. Differential Equations} \textbf{143} (1998), 267--292.

\bibitem{GT01}
\newblock D. Gilbarg, N.S. Trudinger, 
\newblock {``Elliptic partial differential equations of second order".} 
\newblock Reprint of the 1998 edition. Classics in Mathematics, Springer-Verlag, Berlin,  2001.

\bibitem{Ha81}  
\newblock E.I.~Hanzawa,
\newblock \emph{Classical solutions of the Stefan problem},
\newblock {T\^ohoku Math. Jour.}, \textbf{33} (1981), 297--335.

\bibitem{K08} 
\newblock M. Kimura,
\newblock \emph{Geometry of hypersurfaces and moving hypersurfaces in $\R^m$ for the study of moving boundary problems},
\newblock topics in mathematical modeling, 
\newblock J. Necas Center for Mathematical Modeling, Lecture Notes, 4, Prague (2008), 39-93. 

\bibitem{KPW09} 
\newblock M.~K\"ohne, J.~Pr\"uss and M.~Wilke,
\newblock \emph{On quasilinear parabolic evolution equations in weighted $L_p$-spaces},
\newblock {J.Evol.Eqns.}, \textbf{10} (2010), 443-463.

\bibitem{KPW10} 
\newblock M.~K\"ohne, J.~Pr\"uss and M.~Wilke,
\newblock \emph{Qualitative behaviour of solutions for the two-phase Navier-Stokes equations with surface tension},
\newblock {Math.~Ann.}, to appear, 
\newblock {arXiv:1005.1023}.

\bibitem{K02} 
\newblock W. K\"uhnel, 
\newblock {``Differential geometry. Curves-surfaces-manifolds,"}
\newblock  Student Mathematical Library, 16, American Mathematical Society, Providence, RI, 2002.

\bibitem{PSSS11} 
\newblock J.~Pr\"uss, Y.~Shibata, S.~Shimizu and G.~Simonett,
\newblock \emph{On well-posedness of incompressible two-phase flows with phase transition: The case of equal densities},
\newblock {Evol.Eqns.\& Control Th.} \textbf{1} (2012), 171--194.  

\bibitem{PSW11}
\newblock J.~Pr\"uss, G.~Simonett and M.~Wilke,
\newblock \emph{On thermodynamically consistent Stefan problems with variable surface energy}, submitted
\newblock {arXiv:1109.4542}.

\bibitem{PSZ10} 
\newblock J. Pr\"uss, G. Simonett and R. Zacher,
\newblock \emph{Qualitative behaviour of solutions for thermodynamically consistent Stefan problems with surface tension},
\newblock {Arch. Ration. Mech. Anal.} (2012) (DOI) 10.1007/s00205-012-0571-y. 

\end{thebibliography}
\end{document}